(13/2/2024, 2h14)

# Friable averages of complex arithmetic functions

R. de la Bretèche & G. Tenenbaum

**Abstract.** We evaluate friable averages of arithmetic functions whose Dirichlet series is analytically close to some complex power of the Riemann zeta function. We obtain asymptotic expansions resembling those provided by the Selberg-Delange method in the non-friable case. Some application are provided to the friable distribution of the additive function counting the total number of prime factors.

## 1. Introduction and statement of results

Let $P^+(n)$ denote the largest prime factor of a natural integer $n > 1$ and let us agree that $P^+(1) = 1$. We designate by $S(x, y)$ the set of $y$-friable integers not exceeding $x$, i.e.

$$S(x, y) := \{n \leqslant x : P^+(n) \leqslant y\}.$$

We aim here at evaluating friable sums

$$\Psi(x, y; f) := \sum_{n\in S(x,y)} f(n)$$

for some complex arithmetical functions $f$.

To describe the set of relevant functions $f$, we introduce parameters $\beta > 0$, $\mathfrak{c} \in ]0, 1[$, $\delta > 0$, $z \in \mathbb{C}^*$, and assume that $\beta + \delta < 3/5$. We then consider the class $\mathcal{E}_z(\beta, \mathfrak{c}, \delta)$ of those Dirichlet series $\mathcal{F}(s)$ converging for $\sigma = \Re e\, s > 1$ and represented in this half-plane as

$$\mathcal{F}(s) = \zeta(s)^z \mathcal{B}(s) \tag{1·1}$$

where the series $\mathcal{B}(s) := \sum_{n\geqslant 1} b(n)/n^s$ may be holomorphically continued in the domain

$$\mathcal{D}(\beta, \mathfrak{c}, \delta) := \Big\{s \in \mathbb{C} : \sigma > 1 - \mathfrak{c}/\{1 + \log^+ |\tau|\}^{(1-\delta-\beta)/(\beta+\delta)}\Big\}, \tag{1·2}$$

with $\log^+ t := \max(0, \log t)$ $(t > 0)$, and satisfies

$$\mathcal{B}(s) \ll \{1 + |\tau|\}^{1-\delta} \qquad \big(s \in \mathcal{D}(\beta, \mathfrak{c}, \delta)\big), \tag{1·3}$$

$$\mathcal{B}(s, y) := \sum_{P^+(n)\leqslant y} \frac{b(n)}{n^s} = \mathcal{B}(s) + O\Big(\frac{1}{L_{\beta+\delta}(y)}\Big) \quad \begin{pmatrix} y \geqslant 2,\ \sigma > 1 - \mathfrak{c}/(\log y)^{\beta+\delta}, \\ |\tau| \leqslant L_{\beta+\delta}(y) \end{pmatrix}. \tag{1·4}$$

Here and in the sequel we write

$$L_c(y) := \mathrm{e}^{(\log y)^c} \quad (c > 0,\ y \geqslant 1). \tag{1·5}$$

We then consider, for $\kappa \geqslant |z|$, the subclass $\mathcal{H}(z, \kappa; \beta, \mathfrak{c}, \delta)$ comprising those functions $f$ whose associated Dirichlet series belongs to $\mathcal{E}_z(\beta, \mathfrak{c}, \delta)$ and possesses a majorant series

$$\mathcal{F}^\dagger(s) := \sum_{n\geqslant 1} f^\dagger(n)/n^s \in \mathcal{E}_\kappa(\beta, \mathfrak{c}, \delta).$$

We next write

$$\begin{aligned} &\varepsilon_z := \mathrm{sgn}(\Re e\, z), \quad m := \varepsilon_z \lceil \Re e\, z\rceil, \quad \vartheta := \varepsilon_z m - \Re e\, z \in [0, 1[, \\ &m_z := \lceil \Re e\, z\rceil = \varepsilon_z m, \quad \vartheta_z := \vartheta - i\Im m\, z, \end{aligned} \tag{1·6}$$

and note that $z = m_z - \vartheta_z$.

Our results involve some solutions of delay differential equations which we now describe. First define $g_z:\mathbb{R}\to\mathbb{C}$ as the solution of the delay differential system

$$\begin{cases} vg_z'(v)+zg_z(v-1)=0 & (v>1)\\ g_z(v):=0 & (v<0),\\ g_z(v):=1 & (0\leqslant v\leqslant 1).\end{cases}$$

Thus $g_1$ coincides with Dickman's function $\varrho$, whose Laplace transform satisfies

$$\widehat{\varrho}(s):=\int_{\mathbb{R}}\varrho(v)\mathrm{e}^{-vs}\,\mathrm{d}v=\mathrm{e}^{\gamma+I(-s)},\tag{1·7}$$

where $\gamma$ designates Euler's constant and where we have put $I(w):=\int_0^w(\mathrm{e}^v-1)\,\mathrm{d}v/v$. From the general theory displayed in [14] we have

$$\widehat{g_z}(s):=\int_{\mathbb{R}}g(v)\mathrm{e}^{-sv}\,\mathrm{d}v=s^{z-1}\widehat{\varrho}(s)^z\qquad(s\in\mathbb{C}\smallsetminus\mathbb{R}^-),\tag{1·8}$$

$$g_z(v)=\sum_{0\leqslant j\leqslant J}\frac{c_j}{v^{z+j}\Gamma(1-z-j)}+O\Big(\frac{1}{v^{z+J+1}}\Big)\qquad(v>0,\ J\geqslant 0),\tag{1·9}$$

where $\{c_j\}_{j=0}^{\infty}$ is the sequence of Taylor coefficients of $\widehat{\varrho}(s)^z$ at the origin.

When $\Re e\, z>0$, the expression $\widehat{\varrho}(s)^z$ is the Laplace transform of a function $\varrho_z$ solution of the delay differential system

$$\begin{cases} v\varrho_z'(v)+(1-z)\varrho_z(v)+z\varrho_z(v-1)=0 & (v>1),\\ \varrho_z(v)=v^{z-1}/\Gamma(z) & (0<v\leqslant 1).\end{cases}\tag{1·10}$$

Thus, continued by 0 on $]-\infty,0]$, the function $\varrho_z$ is the order $z$ fractional convolution power of Dickman's function $\varrho=\varrho_1$.

From the definition, we see that $\varrho_z$ is $\mathcal{C}^\infty$ on $\mathbb{R}\smallsetminus\mathbb{N}$ and is $\mathcal{C}^{j-1}$ on $]j-m,\infty[$ provided $j\geqslant m$. When $z=m\in\mathbb{N}^*$, the discontinuities of $\varrho_z$ on $\mathbb{N}^*$ are of the first kind. We may then continue $\varrho_z$ by right-continuity. We put

$$\delta_{z,h,j}=\varrho_z^{(j)}(h)-\varrho_z^{(j)}(h-0)\qquad(z=m\in\mathbb{N}^*,\ 1\leqslant h\leqslant j+1-m).\tag{1·11}$$

Let $\zeta_0(w)$ denote the solution with smallest non-negative imaginary part of the equation $\mathrm{e}^\zeta=1+w\zeta$. Note in particular that $\zeta_0(w)\in\mathbb{R}$ whenever $w\in\mathbb{R}^+$; in this case one traditionally notes $\xi(w):=\zeta_0(w)$. By convention $\xi(w)\neq 0$ for $w>0$, $w\neq 1$, and $\xi(1)=0$. From [20; (III.5.48)] and [14; (2.14)], writing $\mathfrak{t}:=\arg(z)\in]-\pi,\pi]$, and $\xi_z:=\xi(v/|z|)$, we have,

$$\begin{aligned}\xi_z&=\log\Big(\frac{v}{|z|}\log v\Big)+O\Big(\frac{\log_2 v}{\log v}\Big),\\ \zeta_0(v/z)&=\xi_z+\frac{\mathfrak{t}^2}{2\xi_z^2}-\frac{i\xi_z\mathfrak{t}}{\xi_z-1}+O\Big(\frac{1}{\xi_z^3}\Big)\end{aligned}\qquad(v\geqslant 3+|z|).\tag{1·12}$$

The asymptotic behaviour of $\varrho_z(v)$ as $v\to\infty$ has been determined in [14] in terms of $\zeta_0(v/z)$. We have

$$\varrho_z(v)=\frac{\{1+O(1/v)\}\mathrm{e}^{\gamma z-v\zeta_0(v/z)+zI(\zeta_0(v/z))}}{\sqrt{2\pi v\{1-1/\zeta_0(v/z)\}}}\qquad(v\to\infty).\tag{1·13}$$

We next define $\varphi_z$ as the solution of the delay differential system

$$\begin{cases} v\varphi_z'(v)+\vartheta_z\varphi_z(v)+z\varphi_z(v-1)=0 & (v>1),\\ \varphi_z(v)=v^{-\vartheta_z}/\Gamma(1-\vartheta_z) & (0<v\leqslant 1),\\ \varphi_z(v)=0 & (v\leqslant 0).\end{cases}\tag{1·14}$$

By Lemma 3.2 *infra*, with notation (1·6), we have $\varphi_z=\varrho_z^{(m-1)}$ whenever $\Re e\, z>0$. When $z\in\mathbb{Z}$, we have $\vartheta_z=0$, hence the discontinuities of $\varphi_z$ on $\mathbb{N}^*$ are of the first kind and we may continue $\varphi_z$ by right continuity. With this convention, we extend definition (1·11) to $z\in\mathbb{Z}^-$ by substituting $\varphi_z$ to $\varrho_z$ in this case and setting the range for $h$ to $1\leqslant h\leqslant j$.

In Lemma 3.2, it is shown that, in full generality,

$$\widehat{\varphi_z}(s) = s^{m_z-1}\widehat{\varrho}(s)^z \qquad (\Re e\, s > 0).$$

Put

$$M(x;f) := \Psi(x,x;f) = \sum_{n\leqslant x} f(n), \tag{1·15}$$

$$\lambda_{y,f}(u) := \int_{\mathbb{R}} g_z(u-v)\,\mathrm{d}\Big(\frac{M(y^v;f)}{y^v}\Big), \quad \Lambda_f(x,y) := x\lambda_{y,f}(u), \tag{1·16}$$

$$R_z(v) := \frac{1}{\sqrt{v}}\exp\Big\{-\Re\int_r^u \zeta_0(t/z)\,\mathrm{d}t\Big\}, \tag{1·17}$$

and let us note incidentally that if, following [14], we let $\{F_k(v;a,b)\}_{k=0}^{\infty}$ denote the sequence of fundamental solutions to the general delay differential equation

$$v\mathfrak{f}'(v) + a\mathfrak{f}(v) + b\mathfrak{f}(v-1) = 0,$$

then we actually have $R_z(v) \asymp F_0(v;-1-z,z)$ $(v \geqslant 1)$.

Our first theorem states that, if $f \in \mathcal{H}(z,\kappa;\beta,\mathfrak{c},\delta)$, then $\Lambda_f(x,y)$ is a good approximation to $\Psi(x,y;f)$ in the range

$$x \geqslant 3, \quad \exp\{(\log x)^{1-\beta}\} \leqslant y \leqslant x. \tag{$G_\beta$}$$

Here and throughout we write $u := (\log x)/\log y$ and systematically employ the notation

$$\begin{aligned} &\sigma := \Re e\, s, \quad \tau := \Im m\, s \quad (s\in\mathbb{C}),\\ &H^+ := \mathbb{R}^{+*} + i\mathbb{R}, \quad H^- := \mathbb{R}^- + i\mathbb{R}. \end{aligned} \tag{1·18}$$

**Theorem 1.1.** *Let $\beta > 0$, $\mathfrak{c} > 0$, $\delta > 0$, $z \in \mathbb{C}^*$, $\kappa > 0$, $\beta + \delta < 3/5$. Then, uniformly for $f \in \mathcal{H}(z,\kappa;\beta,\mathfrak{c},\delta)$ and $(x,y) \in G_\beta$, we have*

$$\Psi(x,y;f) = \Lambda_f(x,y) + O\Big(\frac{xR_z(u)}{L_{\beta+\delta/2}(y)}\Big). \tag{1·19}$$

Next, it is necessary to evaluate $\Lambda_f(x,y)$. The statement of our results addressing this question necessitates a number of further definitions.

We denote by $\mathrm{d}\mu_{y,f}$ the real measure with Laplace transform

$$\widehat{\mu_{y,f}}(s) := \int_{\mathbb{R}} \mathrm{e}^{-vs}\,\mathrm{d}\mu_{y,f}(v) = \mathcal{F}\Big(1+\frac{s}{\log y}\Big)\frac{s^{1-\vartheta_z}}{s+\log y} \qquad (\sigma > 0).$$

We prove in Lemma 3.4 *infra* that $\mathrm{d}\mu_{y,f}(v)$ is absolutely continuous whenever $\Re e\, z \in \mathbb{R}\smallsetminus\mathbb{Z}$ and that, in this case,

$$\mathrm{d}\mu_{y,f}(v) = Z_{y,f}(v)\,\mathrm{d}v, \quad \text{with} \quad Z_{y,f}(v) := \frac{1}{2\pi i}\int_{1+i\mathbb{R}} \widehat{\mu_{y,f}}(s)\mathrm{e}^{vs}\,\mathrm{d}s \quad (y^v \in \mathbb{R}\smallsetminus\mathbb{N}). \tag{1·20}$$

The function $Z_{\mathrm{e},f}$ is made explicit in Lemma 3.4 and we note that

$$Z_{y,f}(v) = (\log y)^{1-\vartheta_z} Z_{\mathrm{e},f}(v\log y).$$

From definition (1·20) we have

$$\widehat{Z_{y,f}}(s) = \frac{\mathcal{F}(1+s/\log y)s^{1-\vartheta_z}}{s+\log y} = \frac{\widehat{Z_{\mathrm{e},f}}(s/\log y)}{(\log y)^{\vartheta_z}} \qquad (\Re e\, s > 0).$$

Thus

$$\begin{aligned} \widehat{\lambda_{y,f}}(s) &= \widehat{g_z}(s)\int_0^\infty \mathrm{e}^{-vs}\,\mathrm{d}\Big(\frac{M(y^v;f)}{y^v}\Big)\\ &= \widehat{g_z}(s)\frac{s\mathcal{F}(1+s/\log y)}{s+\log y} = \widehat{\mu_{y,f}}(s)s^{\vartheta_z}\widehat{g_z}(s), \end{aligned} \tag{1·21}$$

whence

$$\widehat{\lambda_{y,f}}(s) = \widehat{\mu_{y,f}}(s)s^{m_z-1}\widehat{\varrho}(s)^z = \widehat{\mu_{y,f}}(s)\widehat{\varphi_z}(s). \tag{1·22}$$

Let the sequence $\{a_j(f)\}_{j=0}^{\infty}$ be defined by the Taylor expansion converging in a neighbourhood of the origin

$$\frac{s^z\mathcal{F}(s+1)}{s+1}=\sum_{j\geqslant 0}a_j(f)s^j\qquad(|s|<\mathfrak{c})\tag{1·23}$$

and let us make the convention that $a_j(f)=0$ if $j\in\mathbb{Z}^-$. We then define

$$\mathrm{d}\nu_{y,f}(v):=\mathrm{d}\mu_{y,f}(v)-(\log y)^{z-1}\sum_{1\leqslant j<m_z}\frac{a_{m-1-j}(f)v^{j-1}\,\mathrm{d}v}{(j-1)!(\log y)^{m-1-j}}\qquad(v>0),\tag{1·24}$$

so that $\mathrm{d}\nu_{y,f}(v)=\mathrm{d}\mu_{y,f}(v)$ provided $\Re e\, z<0$, and, in complete generality,

$$\int_0^{\infty}\mathrm{e}^{-sv}\,\mathrm{d}\nu_{y,f}(v)=(\log y)^{-\vartheta_z}\sum_{j\geqslant m_z-1}a_j(f)\Big(\frac{s}{\log y}\Big)^{j+1-m_z}\qquad(|s|<\mathfrak{c}\log y).\tag{1·25}$$

Writing

$$W_j(t,y;f):=\int_t^{\infty}(v-t)^j\,\mathrm{d}\nu_{y,f}(v)\qquad(t\geqslant 0,\ j\geqslant 0,\ y\geqslant 1),\tag{1·26}$$

it follows by partial integration from the case $j=0$ of (3·20) *infra* that

$$W_j(t,y;f)\ll\frac{1}{(\log y)^{j+\kappa_z-\vartheta_z}L_b(y^t)}\cdot\tag{1·27}$$

For integer $J\geqslant 0$, real $y\geqslant 3$, we put $e_y:=(\log_2 y)^{1/\beta}/\log y$ and consider the sets

$$\mathcal{V}_J(y):=\Big\{u\geqslant 1:\min_{1\leqslant j\leqslant\min(u,J+1)}(u-j)\geqslant e_y\Big\}.$$

Defining

$$\psi_z:=\begin{cases}\varphi_z^{(m+1)} & \text{if } \Re e\, z\leqslant 0,\\ \varrho_z & \text{if } \Re e\, z>0,\end{cases}\tag{1·28}$$

we may state our main result below in a unified setting.

**Theorem 1.2.** *Let*

$$\beta>0,\ \mathfrak{c}>0,\ \delta>0,\ \beta+\delta<3/5,\ \kappa>0,\ z\in\mathbb{C}^*,$$
$$J\geqslant m_z,\ J_z:=J+1-m_z.$$

*Then, uniformly for $f\in\mathcal{H}(z,\kappa;\beta,\mathfrak{c},\delta)$, $(x,y)\in G_\beta$, $u=(\log x)/\log y\in\mathcal{V}_{J_z}(y)$, we have*

$$\Psi(x,y;f)=x(\log y)^{z-1}\sum_{0\leqslant j\leqslant J}\frac{a_j(f)\psi_z^{(j)}(u)}{(\log y)^j}+O\Big(\frac{xR_z(u)(\log 2u)^{J+1}}{(\log y)^{J+2-z}}\Big).\tag{1·29}$$

*When $u\notin\mathcal{V}_{J_z}(y)$ and $z\in\mathbb{Z}$, the above formula persists provided the quantity $xU_J(x,y;f)$ is added to the main term, with*

$$U_J(x,y;f):=\sum_{\ell\leqslant j\leqslant J_z}\frac{(-1)^{j+1}\delta_{z,\ell,j}}{j!}W_j(u-\ell,y;f)\qquad(\ell<u\leqslant\ell+1).\tag{1·30}$$

*If $u\notin\mathcal{V}_{J_z}(y)$, $\ell<u\leqslant\ell+1$, $\ell\leqslant J_z+1$, and $z\in\mathbb{C}^*\smallsetminus\mathbb{Z}$, formula (1·29) must be modified by restricting the summation to the (possibly empty) range $0\leqslant j<\ell+m_z$ and replacing the error term by $x/(\log y)^{\ell}$.*

*Remarks.* (i) If $z=m\in\mathbb{N}^*$, $m\geqslant 2$, then $\delta_{z,\ell,j}=0$ whenever $\ell\geqslant j-m+2$.

(ii) An inspection of the proofs shows that, given any $\varepsilon>0$, our estimates are uniform for $z$ running in any bounded subset of $H^+\smallsetminus\{|z|<\varepsilon\}$. This is important for the applications described in § 2. For $z\in H^-$, Corollaries 1.3(i) and 1.4 will suffice.

The next corollary is obtained by selecting $J := \max(0, m_z)$ in the previous statement. It extends [5; cor. 1.3] to $z \in H^-$ and [9; (1.24)] to $z \in H^+$. We let $\langle t\rangle$ denote the fractional part of a real number $t$.

**Corollary 1.3.** (i) *Let $z \in H^-$. Then, uniformly for $f \in \mathcal{H}(z,\kappa;\beta,\mathfrak{c},\delta)$, $(x,y) \in G_\beta$, $u = (\log x)/\log y \in \mathcal{V}_{m+1}(y)$, we have*

$$\Psi(x,y;f) = x(\log y)^{z-1}\Big\{a_0(f)\psi_z(u) + O\Big(\frac{xR_z(u)\log 2u}{\log y}\Big)\Big\}. \tag{1·31}$$

*Moreover, if $y^\ell < x \leqslant y^{\ell+1}$ for some integer $\ell \in [1, m+1]$ and $\langle u\rangle \leqslant e_y$, we have*

$$\Psi(x,y;f) \ll \frac{x}{(\log y)^{\min(1-\Re e\, z,\ell)}}. \tag{1·32}$$

(ii) *Let $z \in H^+$. Then, uniformly for $f \in \mathcal{H}(z,\kappa;\beta,\mathfrak{c},\delta)$, $(x,y) \in G_\beta$, we have*

$$\Psi(x,y;f) = x(\log y)^{z-1}\Big\{a_0(f)\psi_z(u) + O\Big(\frac{R_z(u)\log 2u}{\log y} + \frac{\mathbf{1}_{[1,1+e_y]}(u)}{(\log y)^{\Re e\, z}}\Big)\Big\}. \tag{1·33}$$

Combined with some estimates from the literature, our results open the way to upper bounds that are uniform for $x \geqslant y \geqslant 2$. The following corollary generalises [5; cor. 1.4]. It is proved in §6. We denote by $\mathcal{H}^*(z,r;\beta,\mathfrak{c},\delta)$ the subclass of $\mathcal{H}(z,r;\beta,\mathfrak{c},\delta)$ subject to the further condition that, for a suitable constant $C$, we have

$$\mathcal{B}(s,y) \ll \zeta(2\alpha_r,y)^C, \quad 1/\zeta(2\alpha_r,y)^C \ll \mathcal{B}^\dagger(s,y) \ll \zeta(2\alpha_r,y)^C \qquad (\Re e\, s \geqslant \alpha_r), \tag{1·34}$$

where $\alpha_r = \alpha_r(x,y)$ is defined as the saddle-point associated to $\zeta(s,y)^r$.

**Corollary 1.4.** *Let $\beta > 0$, $\mathfrak{c} > 0$, $\delta > 0$, $\beta + \delta < 3/5$, $1 < b < 3/2$, $z \in H^-$, with $|z| = r$. There exists a constant $c_0$ such that, uniformly for $f \in \mathcal{H}^*(z,r;\beta,\mathfrak{c},\delta)$, we have*

$$\Psi(x,y;f) \ll \Psi(x,y;f^\dagger)\Big\{\frac{\mathrm{e}^{-c_0u/(\log 2u)^2}}{(\log y)^{r+d-1}} + \frac{1}{L_b(y)}\Big\} \quad (x \geqslant y \geqslant 2), \tag{1·35}$$

*where $d := \min(\lfloor u\rfloor, 1 - \Re e\, z)$. This bound is locally uniform for $z \in H^- \smallsetminus \{0\}$.*

## 2. Some applications

Our results apply in a natural way to the distribution of additive functions over $S(x,y)$. By way of illustration we state three corollaries, proved in section 7, relevant to the function

$$\omega(n) := \sum_{p\mid n} 1 \qquad (n \geqslant 1).$$

For $r > 0$, we define

$$\begin{aligned} \mu_r = \mu_r(x,y) &:= r\log_2 y + rI\big(\xi(u/r)\big), \\ \sigma_r^2 = \sigma_r(x,y)^2 &:= \mu_r - u^2\xi'(u/r)/r = \mu_r - \frac{ru^2\xi(u/r)}{r + \xi(u/r)(u-r)}, \end{aligned} \tag{2·1}$$

and note that, by (3·7) and (7·6) *infra*, we have, for bounded $r$,

$$\sigma_r^2 = r\log_2 y + \frac{ru}{\xi(u)^2} + O\Big(\frac{u}{\xi(u)^3}\Big) \qquad (u \geqslant r+2).$$

When $r = 1$, we simply write $\mu := \mu_1$, $\sigma := \sigma_1$.

We start with an Erdős-Kac type theorem for friable integers. We use the traditional notation

$$\Phi(v) := \frac{1}{\sqrt{2\pi}}\int_{-\infty}^{v} \mathrm{e}^{-t^2/2}\,\mathrm{d}t \qquad (v \in \mathbb{R})$$

for the distribution function of the normalised Gaussian law.

**Corollary 2.1.** *Let $K>1$. Uniformly for $v\in\mathbb{R}$, $(x,y)\in G_\beta$, we have*

$$\frac{1}{\Psi(x,y)}\sum_{\substack{n\in S(x,y)\\ \omega(n)-\mu\leqslant v\sigma}}1=\Phi(v)+O\Big(\frac{u}{u^K+\sigma}+\frac{1}{\sigma}\Big). \tag{2·2}$$

*Remarks.* (i) This application only uses first term approximations of our estimates in Theorem 1.2. More sophisticated results could be devised using the full expansions.

(ii) When $u=1$, i.e. $y=x$, the above result recovers the classical Erdős-Kac theorem with optimal accuracy—see, e.g. [20; § III.4.4].

(iii) Writing $w(x):=(\log_3 x)^2\log_2 x$, $\gamma(K):=1/(2K+2)$, the remainder term of (2·2) is

$$\ll\begin{cases}\dfrac{\log u}{\sqrt{u}} & \text{if } \mathrm{e}^{(\log x)^{1-\beta}}\leqslant y\leqslant x^{1/w(x)},\\[2mm] \dfrac{\log u}{u^K}+\dfrac{1}{\sqrt{\log_2 y}} & \text{if } x^{1/w(x)}<y\leqslant x^{1/w(x)^{\gamma(K)}},\\[2mm] \dfrac{u}{\sqrt{\log_2 y}} & \text{if } x^{1/w(x)^{\gamma(K)}}<y\leqslant x\ .\end{cases}$$

(iv) This problem has been considered before. Hildebrand [12] established Gaussian convergence with explicit remainder in the range $\log x\geqslant(\log y)^{21}$, which does not intersect $G_\beta$. Alladi [2] tackled the case of larger values of $y$, namely $\exp\big\{(\log_2 x)^{1/\beta}\big\}\leqslant y\leqslant x$, by the method of moments. He showed that (2·2) holds with a remainder $o(1)$, without providing an effective bound.(1) Hensley [10] considered the closely related function $\Omega(n):=\sum_{p^\nu\|n}\nu$ and, in the range $(\log x)^2\leqslant y\leqslant\exp\{(\log x)^{1/4}\}$, obtained uniform estimates for the local laws that yield normal convergence. Finally, Mehdizadeh [15], relying on local behaviour properties of $\Psi(x,y)$, obtained a similar estimate for the range $x^{c(x)/\log_2 x}\leqslant y\leqslant x$, provided $c(x)\to\infty$.

(v) All our results in this section could be adapted to handle the function $\Omega(n)$ with straightforward modifications. We leave this to the reader.

Next, we consider the local laws of the friable distribution of $\omega(n)$, i.e. we aim at evaluating

$$\psi_k(x,y):=\frac{1}{\Psi(x,y)}\sum_{\substack{n\in S(x,y)\\ \omega(n)=k}}1. \tag{2·3}$$

Our next statement splits into two parts.

In the first part, we provide a formula that is uniform for $k\geqslant1$ but only yields an asymptotic formula when $|k-\mu|\leqslant\sigma\sqrt{c\log\sigma}$ with $c<2$.

The second part furnishes an asymptotic formula provided we can select $r$ such that $\mu_r=k$. Formula (3·7) below furnishes an acceptable range for $k$, made explicit in the statement. When this condition is met, we write

$$\begin{aligned}L&:=\log_2 y+I\big(\xi(u/r)\big)\asymp\mu,\\ \mathcal{K}_r(u)&:=\prod_p\Big(1-\frac{1}{p}\Big)^{r-1}\Big(1+\frac{r-1}{p}\Big)\frac{\varrho_r(u)\mathrm{e}^{(1-r)I(\xi(u/r))}\sqrt{\mu_r}}{\varrho(u)\sigma_r}.\end{aligned} \tag{2·4}$$

**Corollary 2.2.** *Let $K>0$.*

(a) *Uniformly for $(x,y)\in G_\beta$, $k\geqslant1$, we have*

$$\psi_k(x,y)=\frac{\mathrm{e}^{-\frac12(\mu-k)^2/\sigma^2}}{\sqrt{2\pi}\sigma}+O\Big(\frac{u}{\sigma^2+u^K\sigma}+\frac{1}{\sigma^2}\Big). \tag{2·5}$$

(b) *Let $c_1$, $c_2$, $c_3$ denote arbitrary positive constants. Uniformly for*

$$(x,y)\in G_\beta,\quad c_1\sigma^2-c_2u/(\log 2u)^2\leqslant k-I\big(\xi(u)\big)\leqslant c_3\sigma^2, \tag{2·6}$$

*and if $r$ is given by the equation $\mu_r=k$, we have, with notation (2·4),*

$$\psi_k(x,y)=\mathcal{K}_r(u)\frac{\mathrm{e}^{-L}L^k}{k!}\Big\{1+O\Big(\frac{u}{\sigma+u^K}+\frac{1}{\sigma}\Big)\Big\}. \tag{2·7}$$

1. Alladi's estimate is stated with different expressions for the mean and variance. However, it can be shown that they differ from ours by an amount $\ll\log_2 3u$, which ensures compatibility between the two statements.

*Remarks.* (i) In the intersection of their respective validity ranges, the remainder term of (2·5), which is $\ll 1/\sigma_{x,y}^c$ for any $c < 2$, is roughly the square of that stated in [10; th. 2(b)].

(ii) By suitable summation over $k$, formula (2·5) plainly provides an Erdős-Kac type theorem. However this turns out to be slightly less precise than (2·2), for which we hence gave a separate proof.

(iii) The range for $k$ covered by Corollary 2.2 includes values as large as powers of $\log y$.

Finally, we propose an estimate for large deviations.

**Corollary 2.3.** *Uniformly for $(x, y) \in G_\beta$, $v \geqslant 0$, and suitable absolute constant $c > 0$, we have*

$$\frac{1}{\Psi(x,y)} \sum_{\substack{n \in S(x,y) \\ |\omega(n)-\mu|>v\sigma}} 1 \ll \mathrm{e}^{-v^2/3} + \mathrm{e}^{-c\sigma^2/(\log\sigma)^4}. \tag{2·8}$$

*Moreover, if $v \ll \min\big(\sigma/u^{1/3}, \sigma^{1/3}\big)$ the upper bound may be replaced by $\ll \mathrm{e}^{-v^2/2}$.*

In the intersection of the respective ranges of validity of the two formulae, the bound (2·8) is significantly more precise than that given in [10; th. 2(a)].

## 3. Lemmas

### *3·1. Bounds of Laplace transforms*

For $z = r\mathrm{e}^{it} \in \mathbb{C}^*$, $\mathrm{t} \in ]-\pi, \pi]$, and $v > 0$, let us put

$$\zeta_z = \zeta_z(v) := \zeta_0(v/z) = \xi_0(v/z) + i\eta_0(v/z), \tag{3·1}$$

$$\xi_z = \xi_z(v) := \xi(v/r). \tag{3·2}$$

**Lemma 3.1.** *For $v \geqslant 1$, $s = -\zeta_z(v) + i\tau$, $\tau \in \mathbb{R}$, we have*

$$\big|s^z \widehat{\varrho}(s)^z\big| \ll \big|\zeta_z(v)^z \widehat{\varrho}(-\zeta_z(v))^z\big|. \tag{3·3}$$

*Proof.* Recall (1·7) and (1·12). The upper bound (3·3) is equivalent to

$$\Re e\,\big\{zI(-s) + z\log s - zI(\zeta_z) - z\log\xi_0(v/z)\big\} \leqslant O(1). \tag{3·4}$$

Since, by [20; lemma III.5.9],

$$I(-s) + \log s = -\gamma - J(s), \quad \text{with} \quad J(s) := \int_0^\infty \frac{\mathrm{e}^{-s-t}}{s+t}\,\mathrm{d}t, \tag{3·5}$$

we may assume that $v$ is arbitrarily large. Indeed, if $v$ is bounded and $|\tau|$ is large, it follows from (3·5) that $I(-s) + \log s$ is bounded and estimates (1·12) imply that $I(\zeta_z)$ and $\log\xi_0(v/z)$ are also bounded.

Writing $T(s) := 1 + 1/s + 2/s^2$, we have by [14; lemma 2] that

$$I(s) = \frac{\mathrm{e}^s}{s}T(s) + O\Big(\frac{\mathrm{e}^\sigma}{\sigma^4}\Big) \qquad (\sigma \geqslant 1), \tag{3·6}$$

whence

$$I(\zeta_z) = \frac{v}{z}T(\zeta_z) + O\Big(\frac{v}{\xi_z^3}\Big) \qquad (v \geqslant 2r). \tag{3·7}$$

By [20; lemma III.5.12], we have $|s\widehat{\varrho}(s)| \asymp 1$ when $|\tau| \geqslant C\{1 + v\xi(v)\}$ with a sufficiently large constant $C$. It hence follows from (3·7) that (3·4) holds in this range. In the complementary circumstance, the term $\log s$ in (3·4) is $\ll \log v$. However, by formulae (1·12) we have

$$|z\mathrm{e}^s| \leqslant |z|\mathrm{e}^{\xi_z}\Big\{1 + O\Big(\frac{1}{\xi_z^2}\Big)\Big\} \leqslant v\xi_z + O\Big(\frac{v}{\xi_z}\Big).$$

As a consequence, in view of (3·6) and (3·7), we see that (3·4) is satisfied if $|\tau|$ is sufficiently large. We may assume in the sequel that $\tau$ is bounded.

Thus, it is enough to establish (3·4) for large $v$ and $\tau \ll 1$. Now

$$\begin{aligned} I(\zeta_z - i\tau) &= \frac{\mathrm{e}^{\zeta_z - i\tau}}{\zeta_z}\frac{\zeta_z}{\zeta_z - i\tau}T(\zeta_z - i\tau) + O\Big(\frac{v}{\xi_z^3}\Big) \\ &= \frac{\mathrm{e}^{-i\tau}v}{z}\frac{\zeta_z}{\zeta_z - i\tau}T(\zeta_z - i\tau) + O\Big(\frac{v}{\xi_z^3}\Big) \\ &= \frac{\mathrm{e}^{-i\tau}v}{z}\Big(1 + \frac{1+i\tau}{\zeta_z} + \frac{2+2i\tau-\tau^2}{\zeta_z^2}\Big) + O\Big(\frac{v}{\xi_z^3}\Big), \end{aligned} \tag{3·8}$$

hence

$$\begin{aligned} &zI(\zeta_z - i\tau) - zI(\zeta_z) \\ &\qquad = v\Big\{\mathrm{e}^{-i\tau} - 1 + \frac{\mathrm{e}^{-i\tau}(1+i\tau)-1}{\zeta_z} + \frac{\mathrm{e}^{-i\tau}(2+2i\tau-\tau^2)-2}{\zeta_z^2} + O\Big(\frac{1}{\xi_z^3}\Big)\Big\}. \end{aligned} \tag{3·9}$$

Gathering estimates (3·6), (3·7) and (3·8), we infer that

$$\begin{aligned} zI(\zeta_z - i\tau) - zI(\zeta_z) &= v\Big\{\mathrm{e}^{-i\tau}\Big(1 + \frac{1+i\tau}{\zeta_z} + \frac{2+2i\tau-\tau^2}{\zeta_z^2}\Big) - T(\zeta_z) + O\Big(\frac{1}{\xi_z^3}\Big)\Big\} \\ &= v\Big\{\mathrm{e}^{-i\tau} - 1 + \frac{\mathrm{e}^{-i\tau}(1+i\tau)-1}{\zeta_z} + \frac{\mathrm{e}^{-i\tau}(2+2i\tau-\tau^2)-2}{\zeta_z^2} + O\Big(\frac{1}{\xi_z^3}\Big)\Big\}. \end{aligned} \tag{3·10}$$

By (1·12), we may replace $\zeta_z$ by $\xi_0(v/z) - i\mathfrak{t}\xi_z/(\xi_z - 1)$ and then by $\xi_z - i\mathfrak{t}$. Indeed,

$$\frac{1}{\zeta_z} = \frac{1}{\xi_z - i\mathfrak{t}} + O\Big(\frac{1}{\xi_z^3}\Big) = \frac{1}{\xi_z} + \frac{i\mathfrak{t}}{\xi_z^2} + O\Big(\frac{1}{\xi_z^3}\Big).$$

The quantity inside curly brackets in the right-hand side of (3·10) is

$$\mathrm{e}^{-i\tau} - 1 + \frac{\mathrm{e}^{-i\tau}(1+i\tau)-1}{\xi_z} + \frac{\mathrm{e}^{-i\tau}(2+2i\tau-\tau^2) - 2 + i\mathfrak{t}(\mathrm{e}^{-i\tau}(1+i\tau)-1)}{\xi_z^2} + O\Big(\frac{1}{\xi_z^3}\Big)$$

We now show that

$$D(\tau,v) := v\Big\{\cos(\tau) - 1 + \frac{A(\tau)}{\xi_z} + \frac{B(\tau)}{\xi_z^2} + O\Big(\frac{1}{\xi_z^3}\Big)\Big\} \leqslant O(1)$$

with $A(\tau) := \cos(\tau) - 1 + \tau\sin(\tau)$, and

$$\begin{aligned} B(\tau) &:= \Re e\,\big(\mathrm{e}^{-i\tau}(2+2i\tau-\tau^2) - 2 + i\mathfrak{t}(\mathrm{e}^{-i\tau}(1+i\tau)-1)\big) \\ &= (2-\tau^2)\cos\tau + 2\tau\sin\tau - 2 + \mathfrak{t}\sin\tau - \mathfrak{t}\tau\cos\tau \\ &= -2(1-\cos\tau) - \tau^2\cos\tau + 2\tau\sin\tau + \mathfrak{t}\sin\tau - \mathfrak{t}\tau\cos\tau. \end{aligned}$$

Write $\tau = 2\pi n + h$, where $n \in \mathbb{Z}$, $|h| \leqslant \pi$ and, since $\tau$ is bounded, $n \ll 1$. If $|h| \geqslant c_0/\sqrt{\xi(v)}$, for a suitable absolute constant $c_0 > 0$, we have $D(\tau,v) \leqslant -2h^2/\pi^2 + O(1/\xi_z) \leqslant 0$. Otherwise we have

$$\frac{D(\tau,v)}{v} = -\frac{h^2}{2} + \frac{-h^2/2 + (2\pi n + h)h}{\xi_z} + O\Big(h^4 + \frac{|h|}{\xi_z} + \frac{1}{\xi_z^2}\Big).$$

If $\pi/\xi_z < |h| \leqslant c_0/\sqrt{\xi(v)}$, then again $D(\tau,v) \leqslant 0$. We may hence write $h = w/\xi_z$ with $|w| \leqslant \pi$. We get

$$D(\tau,v) = vE(w)/\xi_z^2 + O\big(v/\xi_z^3\big)$$

with

$$E(w) := -\tfrac{1}{2}w^2 - 2\pi n w + (-2\pi n + \mathfrak{t})2\pi n = -\tfrac{1}{2}(w + 2\pi n)^2 - 2\pi n(\pi n + \mathfrak{t}).$$

Thus, $E(w) \leqslant -\pi^2/2$ as soon as $n \neq 0$ and (3·4) follows.

We are hence left with the case $n=0$ and $|h|\ll 1/\xi_z$. If $|s-\zeta_z|\ll 1/\xi_z$, we have

$$I(\zeta_z-i\tau)-I(\zeta_z)=-i\tau v/z-\tau^2\int_0^1 I''(\zeta_z-it\tau)(1-t)\,\mathrm{d}t.$$

Observing that $I''(\zeta_z-it\tau)=v/z+O(v/\xi_z)$ for $0\leqslant t\leqslant 1$, we infer that

$$\Re e\,\big(z(I(\zeta_z-i\tau)-I(\zeta_z))\big)=-\tau^2v\big(1+O(1/\xi_z)\big)\leqslant 0$$

provided $v$ is sufficiently large. Since

$$\log|s|-\log\xi_0(v_z)\ll\tau^2/\xi_z^2+1,$$

we see that (3·4) still holds in this last case. □

### 3·2. *Solutions of delay differential equations*

Recall definition (1·14) for $\varphi_z$.

**Lemma 3.2.** *We have*

$$\widehat{\varphi_z}(s)=s^{m_z-1}\widehat{\varrho}(s)^z=s^{\vartheta_z-1}\mathrm{e}^{-zJ(s)}\qquad(\Re e\,s>0).\tag{3·11}$$

*Proof.* Plainly,

$$\widehat{\varphi_z}(s)=\int_0^\infty\varphi_z(v)\mathrm{e}^{-vs}\,\mathrm{d}v=\frac{1}{s}\int_0^\infty\varphi_z\Big(\frac{v}{s}\Big)\mathrm{e}^{-v}\,\mathrm{d}v,$$

whence

$$\begin{aligned}\frac{\mathrm{d}\{s\widehat{\varphi_z}(s)\}}{\mathrm{d}s}&=-\int_0^\infty\frac{v}{s^2}\varphi_z'\Big(\frac{v}{s}\Big)\mathrm{e}^{-v}\,\mathrm{d}v=\frac{\vartheta_z}{s}\int_0^\infty\varphi_z\Big(\frac{v}{s}\Big)\mathrm{e}^{-v}\,\mathrm{d}v+\frac{z}{s}\int_0^\infty\varphi_z\Big(\frac{v}{s}-1\Big)\mathrm{e}^{-v}\,\mathrm{d}v\\&=\big\{\vartheta_z+z\mathrm{e}^{-s}\big\}\widehat{\varphi_z}(s).\end{aligned}$$

Solving this ordinary differential equation, we get, for some constant $K$,

$$s\widehat{\varphi_z}(s)=Ks^{\vartheta_z}\mathrm{e}^{-zJ(s)}.$$

As $s\to+\infty$, we have

$$\widehat{\varphi_z}(s)\sim\frac{K}{s^{1-\vartheta_z}},\quad\widehat{\varphi_z}(s)\sim\frac{1}{s\Gamma(1-\vartheta_z)}\int_0^s\Big(\frac{v}{s}\Big)^{-\vartheta_z}\mathrm{e}^{-v}\,\mathrm{d}v\sim\frac{1}{s^{1-\vartheta_z}},$$

therefore $K=1$. □

Recall the definition of $R_z(v)$ in (1·17). By [14; (2.16)], we have, with notation (3·1),

$$R_z(u)\asymp\frac{\widehat{\varrho}(-\zeta_z)^z\mathrm{e}^{-u\zeta_z}}{\sqrt{u}}\qquad(u\geqslant 1).\tag{3·12}$$

Under condition $\Re e\,z>0$, we readily deduce from (1·10) by real induction that, given any constant $C>0$, we have $\varrho_z(v)\ll\mathrm{e}^{-Cv}$ $(v\geqslant 1)$. By [14; th. 2], it follows that

$$\varrho_z^{(j)}(v)\ll R_z(v)(\log 2v)^j\qquad(v\geqslant 1,\ j\geqslant 0),\tag{3·13}$$

$$\varrho_z^{(j)}(v-w)\ll R_z(v)(\log 2v)^j\mathrm{e}^{w\xi_z}\qquad(v\geqslant 1,\ 0\leqslant w\leqslant v-\tfrac12,\ 0\leqslant j<m),\tag{3·14}$$

$$\varrho_z^{(m)}(v-w)\ll R_z(v)(\log 2v)^m\mathrm{e}^{w\xi_z}\qquad(v\geqslant 1,\ 0\leqslant w\leqslant v-\tfrac32).\tag{3·15}$$

Recall (1·14). When $\Re e\,z\leqslant 0$, the function $\varphi_z^{(m)}$ satisfies the functional equation

$$v\mathfrak{f}'(v)-z\mathfrak{f}(v)+z\mathfrak{f}(v-1)=0\quad\big(v\in\mathbb{R}\smallsetminus\{0,1,\dots,m\}\big).$$

Since the exceptional solution (as defined in [14]) of this equation is the constant function equal to 1, we deduce from [14; th. 2] that, as $v\to\infty$,

$$\varphi_z^{(m+j)}(v)=\delta_{0j}\mathrm{e}^{\gamma z}+O\Big(R_z(v)(\log 2v)^{j-1}\Big)\qquad(j\geqslant 0,\ v\geqslant m+j+1),\tag{3·16}$$

with Kronecker's notation. It is useful to recall here that $\psi_z:=\varphi_z^{(m+1)}$.

### 3·3. *Estimates of friable Dirichlet series*

For notational concision, we set

$$s_y:=(s-1)\log y\quad(s\in\mathbb{C},\ y\geqslant 1),\qquad\mathcal{L}_\tau:=\log(2+|\tau|)\quad(\tau\in\mathbb{R}).\tag{3·17}$$

The next statement follows from [20; lemma III.5.16]. We define

$$\{(s-1)\zeta(s)\}^z:=\mathrm{e}^{zh(s)}\tag{3·18}$$

where $h(s)$ designates, in any connected zero free region of the Riemann zeta function containing the halfplane $\sigma>1$, the branch of the complex logarithm of $(s-1)\zeta(s)$ that is real when $s$ is real and $>1$.

**Lemma 3.3.** *Let $z \in \mathbb{C}$ and $\varepsilon > 0$. The formula*

$$\zeta(s,y)^z = \{(s-1)\zeta(s)\}^z (\log y)^z \widehat{\varrho}(s_y)^z \Big\{1 + O\Big(\frac{1}{L_{3/5-\varepsilon}(y)} + y^{-1/\mathcal{L}_\tau^{2/3+\varepsilon}}\Big)\Big\} \tag{3·19}$$

*holds uniformly in the domain*

$$y \geqslant y_0(\varepsilon), \quad \sigma \geqslant 1 - \frac{1}{(\log y)^{2/5+\varepsilon} + \mathcal{L}_\tau^{2/3+\varepsilon}}, \quad |\tau| \leqslant L_{3/2-3\varepsilon}(y).$$

### *3·4. Measures*

In (1·20), we asserted that, provided $\Re e\, z$ is not an integer, $\mathrm{d}\mu_{y,f}(v)$ is absolutely continuous with derivative some function $Z_{y,f}(v)$. We now substantiate this claim.

**Lemma 3.4.** *Assume $\Re e\, z \notin \mathbb{Z}$. Then the function defined almost everywhere by*

$$Z_{\mathrm{e},f}(v) = \frac{-\mathrm{e}^{-v}}{\Gamma(\vartheta_z)} \sum_{n \leqslant \mathrm{e}^v} f(n) \int_0^{v - \log n} \frac{\mathrm{e}^w \,\mathrm{d}w}{w^{1-\vartheta_z}} + \frac{1}{\Gamma(\vartheta_z)} \sum_{n \leqslant \mathrm{e}^v} \frac{f(n)}{n(v - \log n)^{1-\vartheta_z}}$$

*satisfies $Z_{\mathrm{e},f}(v)\,\mathrm{d}v = \mathrm{d}\mu_{\mathrm{e},f}(v)$.*

*Remark.* It follows from this statement that $Z_{y,f}$ is locally integrable under the indicated assumption.

*Proof.* For $\Re e\, z \notin \mathbb{Z}$, we have $\vartheta = \Re e\, \vartheta_z > 0$. From the formulae

$$\left.\begin{aligned} &\int_0^\infty \mathrm{e}^{-ws}\,\mathrm{d}\Big(\frac{M(\mathrm{e}^w;f)}{\mathrm{e}^w}\Big) = \frac{\mathcal{F}(s+1)}{s+1}, \\ &\frac{1}{\Gamma(t)} \int_0^\infty \mathrm{e}^{-ws} w^{t-1}\,\mathrm{d}w = s^{-t} \end{aligned}\right. \qquad (\Re e\,(s) > 0,\ \Re e\,(t) > 0),$$

we deduce, selecting $t = \vartheta_z$, that $s^{-\vartheta_z}\mathcal{F}(s+1)/(s+1)$ is the Laplace transform of

$$v \mapsto \frac{1}{\Gamma(\vartheta_z)} \int_0^v \frac{M(\mathrm{e}^w;f)}{(v-w)^{1-\vartheta_z}\mathrm{e}^w}\,\mathrm{d}w = \frac{\mathrm{e}^{-v}}{\Gamma(\vartheta_z)} \sum_{n \leqslant \mathrm{e}^v} f(n) \int_0^{v-\log n} \frac{\mathrm{e}^w}{w^{1-\vartheta_z}}\,\mathrm{d}w,$$

whose derivative coincides with $Z_{\mathrm{e},f}(v)$. The required conclusion then follows from the classical formula relating the Laplace transform of a function and that of its derivative. □

In the following statement and throughout the sequel of this paper, we denote by $\|t\|$ the distance from a real number $t$ to the set of integers.

**Lemma 3.5.** *Let $\beta > 0$, $\mathfrak{c} > 0$, $\delta > 0$, $\kappa > 0$, $z \in \mathbb{C}^*$, $\beta + \delta < 3/5$, $f \in \mathcal{H}(z,\kappa;\beta,\mathfrak{c},\delta)$. For any fixed integer $j \geqslant 0$, and uniformly for $y \geqslant 2$, $v \geqslant 0$, we have*

$$\int_v^\infty t^j\,\mathrm{d}\nu_{y,f}(t) \ll \frac{1}{(\log y)^{j+\vartheta} L_{\beta+2\delta/3}(y^v)^2}. \tag{3·20}$$

*Under the additional assumption that $\|y^v\| \gg 1$, $\Re e\, z \notin \mathbb{Z}^-$, $0 \leqslant h \leqslant \frac12$, we have*

$$\int_{v-h}^v \mathrm{d}\nu_{y,f}(t) \ll \begin{cases} \dfrac{h^\vartheta}{L_{\beta+2\delta/3}(y^v/h\log y)} & \text{if } \Re e\, z \notin \mathbb{Z}, \\ h & \text{if } \Re e\, z \in \mathbb{N}^*. \end{cases} \tag{3·21}$$

*Remark.* The discontinuity in the expression of the majorant is only apparent: the exponent $\vartheta$ appearing in (3·21) actually arises as $1 - \langle \Re e\, z\rangle$, which coincides with $\vartheta$ when $\Re e\, z \notin \mathbb{Z}$, but is equal to 1 if $\Re e\, z \in \mathbb{Z}$.

*Proof.* The stated result may be derived by following the proofs of [9; lemmas 3.4 & 3.5] when $\Re e\, z > 0$ and of [5; lemma 5.2] when $\Re e\, z \leqslant 0$. We omit the details. □

## 4. Proof of Theorem 1.1

With notation (3·1), let us put

$$\alpha_z := 1 - \xi_0(u/z)/\log y. \tag{4·1}$$

The parameters $\beta$ and $\delta$ being chosen as indicated in the statement, let us put $b := \beta + \delta/2$, $T := u^{2u} L_b(y)^2$. Writing $\mathcal{F}(s,y) := \mathcal{B}(s,y)\zeta(s,y)^z$, Perron's formula (see, e.g., [20; th. II.2.3]) yields

$$\Psi(x,y;f) = \frac{1}{2\pi i}\int_{\alpha_z - iT^2}^{\alpha_z + iT^2} \frac{\mathcal{F}(s,y)x^s}{s}\,\mathrm{d}s + \mathfrak{R}, \tag{4·2}$$

with

$$\mathfrak{R} \ll \sum_{P(n)\leqslant y} \frac{x^{\alpha_z} f^\dagger(n)}{n^{\alpha_z}(1+T^2|\log(x/n)|)} \ll \frac{x^{\alpha_z}\zeta(\alpha_z,y)^\kappa}{T} + \sum_{|n-x|\leqslant x/T} f^\dagger(n). \tag{4·3}$$

By (1·12), we have

$$\xi_0(u/z) = \xi_z + \frac{\mathfrak{t}^2}{2\xi_z^2} + O\Big(\frac{1}{\xi_z^3}\Big).$$

Lemma 3.3, Lemma 3.1 with $\tau = \eta_0(v/z)$, and estimate (3·12) then furnish

$$x^{\alpha_z}\zeta(\alpha_z,y)^\kappa \asymp x(\log y)^\kappa \mathrm{e}^{-u\xi_0(u/z)}\widehat{\varrho}(-\xi_0(u/z))^\kappa \ll x(\log y)^\kappa R_z(u)\mathrm{e}^{O(u)}. \tag{4·4}$$

Morever, since $T \gg x^{1-\varepsilon}$, we have, by Shiu's theorem [16],

$$\sum_{|n-x|\leqslant x/T} f^\dagger(n) \ll x\frac{(\log x)^{\kappa-1}}{T}\cdot \tag{4·5}$$

Carrying these bounds back into (4·3), we get

$$\mathfrak{R} \ll xR_z(u)/L_b(y).$$

The integrand in (4·2) may be estimated appealing to Lemma 3.3, choosing $\beta + \delta < \frac{3}{5} - \varepsilon$. Taking hypothesis (1·4) into account, we get, with notations (3·17),

$$\frac{1}{2\pi i}\int_{\alpha_z - iT^2}^{\alpha_z + iT^2} \mathcal{F}(s,y)\frac{x^s}{s}\,\mathrm{d}s = \frac{(\log y)^z}{2\pi i}\int_{\alpha_z - iT^2}^{\alpha_z + iT^2} \mathcal{F}(s)(s-1)^z \widehat{\varrho}(s_y)^z \frac{x^s}{s}\,\mathrm{d}s + R_1 + R_2,$$

where $\mathcal{F}(s)(s-1)^z$ is defined by means of (3·18) and where the error terms may be estimated appealing to Lemma 3.1 and to estimate (3·12): since $\Re e\,\zeta_z = -\xi_0(u/z) = (\alpha_z - 1)\log y$, we thus get

$$\begin{aligned} R_1 &\ll \frac{x^{\alpha_z}(\log T)^\kappa}{L_{\beta+\delta}(y)} \sup_{\Re e\, s=\alpha_z} |s_y^z \widehat{\varrho}(s_y)^z| \ll \frac{x^{\alpha_z}(\log T)^\kappa}{L_{\beta+\delta}(y)}\big|\zeta_z(u)^z\widehat{\varrho}(-\zeta_z(u))^z\big| \\ &\ll \frac{x\sqrt{u}R_z(u)(\log T)^\kappa}{L_{\beta+\delta}(y)} \ll \frac{xR_z(u)}{L_b(y)}\cdot \end{aligned}$$

and

$$R_2 \ll x^{\alpha_z}(\log T)^\kappa \sup_{\substack{s=\alpha_z+i\tau\\ |\tau|\leqslant T^2}} \Big(|s_y^z\widehat{\varrho}(s_y)^z| y^{-1/\mathcal{L}_\tau^{2/3+\varepsilon}}\Big).$$

To establish these bounds, we made use of the fact that, in a connected zero-free region of $\zeta(s)$ containing the half-plane $\sigma > 1$ and excluding the half-line $]-\infty, 1]$, we have

$$|\{(s-1)\zeta(s)\}^z| = |\zeta(s)^z||(s-1)^z|,$$

where the left-hand side is defined by (3·18) and the right-hand side is defined by means of a holomorphic branch of $\log\zeta(s)$ and the principal determination of $(s-1)^z$.

In order to estimate $R_2$ we use (3·3) when $|\tau| \leqslant y$ and the bound

$$|s_y^z\widehat{\varrho}(s_y)^z| \ll 1 \ll \big|\zeta_z(u)^z\widehat{\varrho}(-\zeta_z(u))^z\big|\mathrm{e}^{-u/2}$$

when $|\tau| > y$. For a suitable choice of $\varepsilon > 0$, we get

$$R_2 \ll xR_z(u)\bigg(\frac{1}{L_b(y)} + \mathrm{e}^{-u/2}y^{-1/\mathcal{L}_{T^2}^{2/3+\varepsilon}}\bigg) \ll \frac{xR_z(u)}{L_b(y)}\cdot$$

At this stage we have hence established the estimate

$$\Psi(x,y;f) = \frac{(\log y)^z}{2\pi i}\int_{\alpha_z-iT^2}^{\alpha_z+iT^2} \mathcal{F}(s)(s-1)^z\widehat{\varrho}(s_y)^z\frac{x^s}{s}\,\mathrm{d}s + O\bigg(\frac{xR_z(u)}{L_b(y)}\bigg). \tag{4·6}$$

The residue theorem allows us to move the integration segment until abscissa $\sigma_x := 1+1/\log x$ at the cost of a satisfactory error term. Indeed, the contribution of the horizontal segments with ordinates $\pm T^2$ is

$$\ll x/T^{2/3} \ll xR_z(u)/L_b(y).$$

Next, we extend the integration segment to the whole vertical line $\Re e\, s = \sigma_x$ taking account of the estimate

$$s\widehat{\varrho}(s) = 1 + O\Big(\frac{u}{1+u|\tau|}\Big) \qquad (\Re e\, s = 1/u) \tag{4·7}$$

and of the bound $|\mathcal{F}(s)| \ll (\log x)^\kappa$ ($\Re e\, s = \sigma_x$). The contribution of the error term in (4·7) is then

$$\ll x(\log x)^\kappa/T \ll xR_z(u)/L_b(y).$$

It remains to bound

$$\frac{1}{2\pi i}\int_{\substack{\sigma=\sigma_x\\|\tau|>T^2}} \mathcal{F}(s)\frac{x^s}{s}$$

using the effective Perron formula as stated, for instance, in [20; th. 2.3]. One checks that this quantity may be absorbed by the previous error term. We have thus proved that

$$\Psi(x,y;f) = \frac{(\log y)^z}{2\pi i}\int_{\sigma_x-i\infty}^{\sigma_x+i\infty} \mathcal{F}(s)(s-1)^z\widehat{\varrho}(s_y)^z\frac{x^s}{s}\,\mathrm{d}s + O\bigg(\frac{xR_z(u)}{L_b(y)}\bigg). \tag{4·8}$$

Since it follows from (1·8) and (1·21) that

$$\frac{(\log y)^z\mathcal{F}(s)(s-1)^z\widehat{\varrho}(s_y)^z x^s}{s} = x\mathrm{e}^{us_y}\frac{s_y\mathcal{F}(1+s_y/\log y)s_y^{z-1}\widehat{\varrho}(s_y)^z}{1+s_y/\log y} = x\mathrm{e}^{us_y}\widehat{\lambda_{y,f}}(s_y),$$

we conclude by the convolution theorem and Laplace inversion that the main term in (4·8) coincides with $\Lambda_f(x,y)$. □

## 5. Proof of Theorem 1.2

### *5·1. Preparation*

In view of Theorem 1.1, it is sufficient to evaluate $\Lambda_f(x,y)$. First observe that it follows from (1·22) that

$$\Lambda_f(x,y) = x\int_{\mathbb{R}} \varphi_z(u-v)\,\mathrm{d}\mu_{y,f}(v). \tag{5·1}$$

Recall the definition of the sequence $\{a_j(f)\}_{j=0}^\infty$ in (1·23). From (1·25), we infer that, still with the convention that $a_h(f) = 0$ for $h < 0$,

$$\frac{a_{j+m_z-1}(f)}{(\log y)^{j+\vartheta_z}} = \frac{(-1)^j}{j!}\int_0^\infty v^j\,\mathrm{d}\nu_{y,f}(v) \qquad (j \geqslant 0). \tag{5·2}$$

Finally, we may assume in all the sequel that $u > 1$ since the Selberg–Delange method as displayed in [20; ch. II.5] provides the required estimates when $u = 1$, i.e. $x = y$.

### 5·2. The case $z \in H^+$

Observe at the outset that we may assume $\|x\| = \frac{1}{2}$.

Apply (5·1) with $\Re e\, z > 0$ and $\varphi_z = \varrho_z^{(m-1)}$. Writing

$$b_j(u) := \int_0^u v^j \varrho_z^{(m-1)}(u-v)\,\mathrm{d}v = j!\varrho_z^{(m-j-2)}(u) \qquad (0 \leqslant j < m-1),$$

we deduce from (1·24) that

$$\frac{\Lambda_f(x,y)}{x(\log y)^{z-1}} = \sum_{1\leqslant j<m} \frac{a_{m-j-1}(f)b_{j-1}(u)}{(j-1)!(\log y)^{m-j-1}} + \mathfrak{J} = \sum_{0\leqslant j<m-1} \frac{a_j(f)\varrho_z^{(j)}(u)}{(\log y)^j} + \mathfrak{J}, \tag{5·3}$$

with

$$\mathfrak{J} := \int_{0-}^u \frac{\varrho_z^{(m-1)}(u-v)}{(\log y)^{z-1}}\,\mathrm{d}\nu_{y,f}(v) = \mathfrak{J}_1 + \mathfrak{J}_2 + \mathfrak{J}_3, \tag{5·4}$$

where the $\mathfrak{J}_\ell$ correspond to the respective integration ranges $[0-, \frac{1}{2}e_y[$, $[\frac{1}{2}e_y, u-\frac{1}{2}[$, $[u-\frac{1}{2}, u]$. Since $\varrho_z^{(m-1)}(w)$ is bounded for $w \geqslant \frac{1}{2}$, the integrals $\mathfrak{J}_1$ and $\mathfrak{J}_2$ are trivially convergent. Moreover, since $\Re e\, z = m - \vartheta$, $0 \leqslant \vartheta < 1$, this is also true for $\mathfrak{J}_3$.

We handle $\mathfrak{J}_2$ and $\mathfrak{J}_3$ as error terms. We retain the notation $b := \beta + \delta/2$.

Let us first consider $\mathfrak{J}_2$. Recall that $\varrho_z^{m-1}$ and that $\varphi_z$ is $\mathfrak{C}^0$ on $]0,\infty[$, and $\mathfrak{C}^1$ on $]1,\infty[$. Write $w := \max(e_y, u-1)$. Partial integration furnishes $\mathfrak{J}_2 = (\mathfrak{J}_{20} + \mathfrak{J}_{21})/(\log y)^{z-1}$, with

$$\mathfrak{J}_{20} := \int_{e_y/2}^w \varphi_z(u-v)\,\mathrm{d}\nu_{y,f}(v) = -\int_{e_y/2}^w \int_{u-v}^\infty \varphi_z'(t)\,\mathrm{d}t\,\mathrm{d}\nu_{y,f}(v), \tag{5·5}$$

$$\mathfrak{J}_{21} := \int_w^{u-1/2} \frac{\mathrm{d}\nu_{y,f}(v)}{\Gamma(1-\vartheta_z)(u-v)^{\vartheta_z}}\cdot \tag{5·6}$$

Inverting the order of integration in (5·5) yields

$$\mathfrak{J}_{20} = -\int_1^\infty \varphi_z'(t)\int_{\max(u-t,e_y/2)}^{\max(u-1,e_y/2)} \mathrm{d}\nu_{y,f}(v)\,\mathrm{d}t \ll \frac{1}{L_b(y^{e_y/2})} \ll \frac{R_z(u)}{(\log y)^{J+2-z}}$$

by (3·14) with $j = m-1$ and (3·20) with $j = 0$.

To evaluate $\mathfrak{J}_{21}$, we apply directly (3·20) with $j = 0$ if $\vartheta_z = 0$. Otherwise we write

$$\begin{aligned}
\mathfrak{J}_{21} &= \int_w^{u-1/2} \Big\{1 + \vartheta_z \int_{u-v}^1 \frac{\mathrm{d}t}{t^{1+\vartheta_z}}\Big\}\frac{\mathrm{d}\nu_{y,f}(v)}{\Gamma(1-\vartheta_z)} \\
&= \int_w^{u-1/2} \frac{\mathrm{d}\nu_{y,f}(v)}{\Gamma(1-\vartheta_z)} + \vartheta_z \int_{1/2}^1 \frac{\mathrm{d}t}{t^{1+\vartheta_z}} \int_{\max(u-t,u-1,e_y/2)}^{u-1/2} \frac{\mathrm{d}\nu_{y,f}(v)}{\Gamma(1-\vartheta_z)} \\
&\ll \frac{1}{L_b(y^{e_y/2})} \ll \frac{R_z(u)}{(\log y)^{J+2-z}},
\end{aligned}$$

still by (3·20) with $j = 0$

We may thus state that we have in all circumstances

$$\mathfrak{J}_2 \ll \frac{R_z(u)}{(\log y)^{J+1}}\cdot \tag{5·7}$$

Let us then turn our attention to $\mathfrak{J}_3$. Assume first that $\vartheta > 0$. Partial integration then yields

$$\begin{aligned}
\mathfrak{J}_3 &= \frac{1}{\Gamma(1-\vartheta_z)(\log y)^{z-1}} \int_{u-1/2}^u \frac{\mathrm{d}\nu_{y,f}(v)}{(u-v)^{\vartheta_z}} \\
&= \frac{1}{\Gamma(1-\vartheta_z)(\log y)^z}\Big\{\Big[-\int_v^u \frac{\mathrm{d}\nu_{y,f}(t)}{(u-v)^{\vartheta_z}}\Big]_{u-1/2}^u + \int_{u-1/2}^u \int_v^u \frac{(\vartheta_z+1)\,\mathrm{d}\nu_{y,f}(t)}{(u-v)^{1+\vartheta_z}}\,\mathrm{d}v\Big\}.
\end{aligned} \tag{5·8}$$

The last integral is estimated using (3·21). It is

$$\ll \int_0^{1/2} \frac{\mathrm{d}h}{hL_b(x/h\log y)} \ll \frac{R_z(u)}{(\log y)^{J+1-z}}.$$

When $\vartheta = 0$, we obtain the required bound from (3·20) in the case $z = m \in \mathbb{N}^*$ since then $\varrho_z^{(m-1)}$ is constant in $]0, 1]$. In the case $z \notin \mathbb{N}^*$, we reach the same conclusion by appealing to the second bound in (3·21), splitting the integral in (5·8) at $u - 1/L_b(y)$.

To evaluate $\mathfrak{J}_1$, we consider two cases according to whether the condition $u \in \mathcal{V}_{J_z}(y)$ is met. If it holds, the function $\varrho_z$ belongs to the class $\mathcal{C}^{J+1}$ on $[u - \frac{1}{2}e_y, u + 0]$. For $0 \leqslant v \leqslant \frac{1}{2}e_y$, the Taylor-Lagrange formula may be written as

$$\varrho_z^{(m-1)}(u-v) = \sum_{m-1\leqslant j\leqslant J} \frac{(-1)^{j-m+1}}{(j-m+1)!}\varrho_z^{(j)}(u)v^{j-m+1} + \mathfrak{R}_0, \tag{5·9}$$

with

$$\mathfrak{R}_0 := \frac{(-1)^{J-m}}{(J-m+1)!}\int_0^v (v-w)^{J-m+1}\varrho_z^{(J+1)}(u-w)\,\mathrm{d}w,$$

while representation (5·2) and estimate (3·20) imply, for $j \geqslant m$,

$$\begin{aligned}\frac{(-1)^{j-m+1}}{(j-m+1)!}\int_{0-}^{e_y/2} v^{j-m+1}\,\mathrm{d}\nu_{y,f}(v) &= \frac{a_j(f)}{(\log y)^{j-m+1+\vartheta_z}} + O\Big(\frac{1}{L_b(y^{e_y/2})}\Big)\\ &= \frac{a_j(f)}{(\log y)^{j+1-z}} + O\Big(\frac{1}{(\log y)^{J+1}}\Big).\end{aligned}\tag{5·10}$$

Carrying back (5·9) into $\mathfrak{J}_1$ taking (5·10) and (3·14) into account, we get

$$\mathfrak{J}_1 = \sum_{m-1\leqslant j\leqslant J} \frac{a_j(f)\varrho_z^{(j)}(u)}{(\log y)^j} + O\bigg(R_z(u)\Big\{\frac{\log(u+1)}{\log y}\Big\}^{J+1}\bigg).$$

It remains to consider the case $u \notin \mathcal{V}_{J_z}(y)$.

If $z = m \in \mathbb{N}^*$, $u > 1$, we have $\ell < u \leqslant \ell + 1$ for a suitable integer $\ell \in [1, J_z]$. The Taylor-Lagrange formula may then be extended by taking account of the first kind discontinuities of $\varrho_z^{(m)}$ and its derivatives. This means adding

$$\mathfrak{R}_1 := \sum_{m\leqslant j\leqslant J+1} \frac{(-1)^{j+1}}{j!} \sum_{\substack{1\leqslant h\leqslant j\\ u-v<h<u}} \delta_{z,h,j}(v+h-u)^j$$

to the remainder $\mathfrak{R}_0$ of (5·9). The quantity

$$\mathfrak{J}_1^\dagger := \frac{1}{(\log y)^{z-1}} \sum_{m\leqslant j\leqslant J+1} \frac{(-1)^{j+1}}{j!} \sum_{\substack{1\leqslant h\leqslant j\\ u-e_y/2<h<u}} \delta_{m,h,j}\int_{u-h}^{e_y/2} (v+h-u)^j\,\mathrm{d}\nu_{y,f}(v)$$

must hence be added to $\mathfrak{J}_1$. The inner sum is actually reduced to the sole term of index $h = \ell$, and estimate (3·20) allows the integration range to be extended to the half-line $[u - \ell, \infty[$ at the cost of an error term compatible with that of (1·29). It follows that

$$(\log y)^{z-1}\mathfrak{J}_1^\dagger = \sum_{m\leqslant j\leqslant J+1} \frac{(-1)^{j+1}}{j!} \sum_{\substack{1\leqslant h\leqslant j\\ u-e_y/2<h<u}} \delta_{m,h,j}W_j(u-h,y;f) + O\Big(\frac{1}{(\log y)^{J+1}}\Big).$$

In view of the previous estimates for $\mathfrak{J}_2$ and $\mathfrak{J}_3$, this completes the proof of (1·29) in the case under consideration.

Let us now turn our attention the case $z \in H^+ \smallsetminus \mathbb{N}^*$, $u \notin \mathcal{V}_{J_z}(y)$. Let then $\ell \in [1, J_z]$ satisfy $\ell < u \leqslant \ell + 1$. Since $\varrho_z^{(m+\ell-1)}$ is integrable on $]\ell, \ell+1]$, we may write the Taylor expansion

$$\varrho_z^{(m-1)}(u-v) = \sum_{m-1 \leqslant j < m+\ell-1} \frac{(-v)^{j-m+1}\varrho_z^{(j)}(u)}{(j-m+1)!} + \mathfrak{R}_0^\dagger \qquad (0 \leqslant v \leqslant \tfrac{1}{2}e_y), \tag{5·11}$$

with

$$\mathfrak{R}_0^\dagger := \frac{(-1)^\ell}{(\ell-1)!} \int_0^v (v-t)^{\ell-1} \varrho_z^{(m+\ell-1)}(u-t)\,\mathrm{d}t.$$

Taking (1·25) and (3·20) into account, we see that the contribution of the sum of (5·11) to the integral $\mathcal{J}_1$ of (5·4) is well approximated by

$$\sum_{m-1 \leqslant j < m+\ell-1} \frac{a_j(f)\varrho_z^{(j)}(u)}{(\log y)^j},$$

the involved error term being compatible with that of (1·29). The contribution of $\mathfrak{R}_0^\dagger$ to $\mathcal{J}_1$ is then

$$\ll \frac{1}{(\log y)^{z-1}} \int_0^{e_y/2} \varrho_z^{(m+\ell-1)}(u-t) \int_t^{e_y/2} (v-t)^{\ell-1}\,\mathrm{d}\nu_{y,f}(v)\,\mathrm{d}t \tag{5·12}$$

We may bound the inner integral by appealing to (3·20). When $\ell > 1$, we have

$$\begin{aligned}\int_t^{e_y/2} (v-t)^{\ell-2}\,\mathrm{d}\nu_{y,f}(v) &= (1-\ell) \int_t^{e_y/2} (v-t)^{\ell-2} \int_v^{e_y/2} \mathrm{d}\nu_{y,f}(w)\,\mathrm{d}v \\ &\ll \int_t^{e_y/2} \frac{(v-t)^{\ell-2}\,\mathrm{d}v}{(\log y)^\vartheta L_{\beta+2\delta/3}(y^v)^2} \ll \frac{1}{(\log y)^{\ell-1+\vartheta} L_b(y^t)^2}.\end{aligned}$$

This estimate also holds for $\ell = 1$. The contribution to the right-hand side of (5·12) from the interval $0 \leqslant t \leqslant \langle u \rangle$ is hence

$$\ll \frac{1}{(\log y)^{z-1}} \int_0^{\langle u \rangle} \frac{(\langle u \rangle - t)^{-\vartheta}\,\mathrm{d}t}{(\log y)^{\ell-1+\vartheta} L_b(y^t)} \ll \frac{1}{(\log y)^{z+\ell-1}}.$$

Indeed, we have

$$\begin{aligned}&\int_0^{\langle u \rangle} \frac{(\langle u \rangle - t)^{-\vartheta}}{L_b(y^t)}\,\mathrm{d}t \\ &\ll \min\Big\{ \int_0^{\langle u \rangle} \frac{\mathrm{d}t}{(\langle u \rangle - t)^\vartheta}, \int_0^{\langle u \rangle/2} \frac{\mathrm{d}t}{\langle u \rangle^\vartheta L_b(y^t)} + \frac{\langle u \rangle^{1-\vartheta}}{L_b(y^{\langle u \rangle/2})} \Big\} \\ &\ll \min\Big\{ \langle u \rangle^{1-\vartheta}, \frac{1}{\langle u \rangle^\vartheta \log y} + \frac{1}{(\log y)^{1-\vartheta}} \Big\} \ll \frac{1}{(\log y)^{1-\vartheta}}.\end{aligned} \tag{5·13}$$

The bound $\varrho_z^{(m+\ell)}(u-t) \ll 1$ may then be employed to estimate he contribution of the interval $\langle u \rangle \leqslant t \leqslant \frac{1}{2}e_y$. It is

$$\ll \frac{1}{(\log y)^{z-1}} \int_{\langle u \rangle}^{e_y/2} \frac{\mathrm{d}t}{(\log y)^{\ell-1+\vartheta} L_b(y^t)^2} \ll \frac{1}{(\log y)^{z+\ell-1+\vartheta} L_b(y^{\langle u \rangle})}.$$

### *5·3. The case $z \in H^-$*

The computations being modelled on those of the previous section, we limit the exposition to brief indications.

Note right away that we may assume $z \notin \mathbb{R}^-$, in view of the results of [5]. Moreover we may assume as before that $\|x\| = \frac{1}{2}$.

In what follows, we restrict to the case $\vartheta > 0$. The case $\vartheta = 0$, $\vartheta_z \neq 0$, may be handled similarly, replacing $\mathrm{d}\mu_{y,f}(v)$ by $\int_v^\infty \mathrm{d}\mu_{y,f}(t)\,\mathrm{d}v$, which has Laplace transform $\widehat{\mu_{y,f}}(s)/s$ and substituting $\varphi_z'$ to $\varphi_z$.

Apply (5·1) with $\Re e\, z \leqslant 0$ and break down

$$\mathfrak{J}^* := \int_0^u \varphi_z(u-v)\,\mathrm{d}\mu_{y,f}(v) = (\log y)^{z-1}\big(\mathfrak{J}_1^* + \mathfrak{J}_2^* + \mathfrak{J}_3^*\big),$$

where the $\mathfrak{J}_\ell^*$ correspond to the respective integration ranges $[0-, \frac12 e_y]$, $]\frac12 e_y, u-\frac12[$, $[u-\frac12, u]$.

The integrals $\mathfrak{J}_2^*$ et $\mathfrak{J}_3^*$ admit upper bounds identical to those leading to (5·7) and (5·8), substituting $\varphi_z$ to $\varrho_z^{(m-1)}$.

To evaluate $\mathfrak{J}_1^*$, we first consider the case $u \in \mathcal{V}_{J_z}(y)$. We may then write Taylor's formula

$$\varphi_z(u-v) = \sum_{0\leqslant j\leqslant J_z} \frac{(-v)^j \varphi_z^{(j)}(u)}{j!} + \mathfrak{R}_0^*,$$

with

$$\mathfrak{R}_0^* := \frac{(-1)^{J_z+1}}{J_z!}\int_0^v (v-t)^{J_z}\varphi_z^{(J_z+1)}(u-t)\,\mathrm{d}t.$$

Appealing to (3·16), (3·20), (5·2), and arguing as in the case $z \in H^+$, we get

$$\begin{aligned}\mathfrak{J}_1^* &= \sum_{0\leqslant j\leqslant J+m+1} \frac{a_{j-m-1}(f)\varphi_z^{(j)}(u)}{(\log y)^{j+\vartheta+z-1}} + O\Big(\frac{R_z(u)(\log 2u)^{J+1}}{(\log y)^{J+1}}\Big)\\ &= \sum_{0\leqslant j\leqslant J} \frac{a_j(f)\psi_z^{(j)}(u)}{(\log y)^j} + O\Big(\frac{R_z(u)(\log 2u)^{J+1}}{(\log y)^{J+1}}\Big).\end{aligned}$$

In the complementary case $u \notin \mathcal{V}_{J_z}(y)$, there exists $\ell \in [1, J_z]$ such that $\ell < u \leqslant \ell+1$. Taylor's formula at order $\ell$ may then be written as

$$\varphi_z(u-v) = \sum_{0\leqslant j<\ell} \frac{(-v)^j\varphi_z^{(j)}(u)}{j!} + \frac{(-1)^\ell}{(\ell-1)!}\int_0^v (v-t)^{\ell-1}\varphi_z^{(\ell)}(u-t)\,\mathrm{d}t.$$

In view of (3·16), (3·20), (5·2), we see that the contribution of the above sum to $\mathfrak{J}_1^*$ is well approximated by

$$\sum_{0\leqslant j\leqslant \ell-m-1} \frac{a_j(f)\psi_z^{(j)}(u)}{(\log y)^j},$$

the summation being empty if $\ell < m+1$, that is to say $u \leqslant m+2$. The contribution of the integral is

$$\mathfrak{R}_1^* := \frac{(-1)^\ell}{(\ell-1)!}\int_0^{e_y/2}\varphi_z^{(\ell)}(u-t)\int_t^{e_y/2}(v-t)^{\ell-1}\frac{\mathrm{d}\nu_{y,f}(v)}{(\log y)^{z-1}}\,\mathrm{d}t.$$

This quantity may be majorised by appealing to the estimate (3·20) and noting that

$$\varphi_z^{(\ell)}(u-t) \ll \begin{cases} (\langle u\rangle - t)^{-\vartheta} & \text{if } 0 \leqslant t \leqslant \langle u\rangle,\\ 1 & \text{if } \langle u\rangle < t \leqslant \frac12 e_y.\end{cases}$$

Indeed we have

$$\begin{aligned}\int_t^{e_y/2}(v-t)^{\ell-1}\,\mathrm{d}\nu_{y,f}(v) &\ll \int_t^{e_y/2}\frac{(v-t)^{\ell-2}\,\mathrm{d}v}{(\log y)^\vartheta L_{\beta+2\delta/3}(y^v)}\\ &\ll \int_t^{e_y/2}\frac{(v-t)^{\ell-2}\,\mathrm{d}v}{(\log y)^\vartheta L_{\beta+\delta/2}(y^t)L_{\beta+\delta/2}(y^{v-t})}\\ &\ll \frac{1}{L_{\beta+\delta/2}(y^t)(\log y)^{\ell-1+\vartheta}}.\end{aligned}$$

By (5·13), it follows that

$$\mathfrak{R}_1^* \ll \frac{1}{(\log y)^{\ell-1+\vartheta+z-1}}\bigg\{\int_0^{\langle u\rangle}\frac{(\langle u\rangle-t)^{-\vartheta}}{L_b(y^t)}\,\mathrm{d}t + \int_{\langle u\rangle}^{e_y/2}\frac{1}{L_{\beta+\delta/2}(y^t)}\,\mathrm{d}t\bigg\} \ll \frac{1}{(\log y)^{\ell+z-1}}.$$

This completes the proof of Theorem 1.2.

## 6. Proof of Corollary 1.4

Throughout this section we denote by $c_j$ $(j = 1, 2, \ldots)$ positive constants depending at most on the parameters $\beta$, $\mathfrak{c}$, $\delta$, $b$, $z$, and ae locally uniform in $r = |z|$.

Let us first assume $(x, y) \in G_\beta$. By [21; th. 2.1] and our assumption (1·34), we then have

$$\Psi(x, y; f^\dagger) \asymp x\varrho_r(u)(\log y)^{r-1},$$

since $\mathcal{B}^\dagger(\alpha_\kappa) \asymp 1$ in this range. Moreover, by (1·13) and (1·12), we have, with $\mathfrak{t} = \arg z$,

$$\frac{R_z(u)}{\varrho_r(u)} \ll \mathrm{e}^{-(r\mathfrak{t}^2/2+o(1))u/(\log 2u)^2}.$$

The required bound (1·35) hence follows from (1·31) and (3·16) if $u \geqslant m+2$, and from (1·32) if $u < m+2$.

Next we turn our attention to the case $u > (\log y)^{\beta/(1-\beta)}$. We then redefine $\alpha_r = \alpha_r(x, y)$ as the saddle point associated to the Perron integral for $\Psi(x, y; \tau_r)$, thereby involving $\zeta(s, y)^r$. A standard application of the saddle-point method as in [6] then furnishes, taking (1·34) into account,

$$\Psi(x, y; f^\dagger) \gg \frac{x^{\alpha_r}\zeta(\alpha_r, y)^r}{\sqrt{\overline{u}}\log y\,\zeta(2\alpha_r, y)^C} \gg \frac{\mathrm{e}^{-c_1\sqrt{\overline{u}}}x^{\alpha_r}\mathcal{F}^\dagger(\alpha_r, y)}{\log y}, \tag{6·1}$$

with the traditional notation $\overline{u} := \min(u, y/\log y)$.

Let $b_1 := \frac{1}{2}(b+\frac{3}{2})$. Evaluating $\Psi(x, y; f)$, as in (4·2), by Perron's formula with $T := L_{b_1}(y)^2$, we get

$$\Psi(x, y; f) = \frac{1}{2\pi i}\int_{\alpha_r - iT^2}^{\alpha_r + iT^2} \frac{\mathcal{F}(s, y)x^s}{s}\,\mathrm{d}s + \mathfrak{R}, \tag{6·2}$$

with

$$\mathfrak{R} \ll \frac{x^{\alpha_r}\mathcal{F}^\dagger(\alpha_r, y)}{T} + \Psi\Big(x + \frac{x}{T}, y; f^\dagger\Big) - \Psi\Big(x - \frac{x}{T}, y; f^\dagger\Big).$$

The first term may be bounded from (6·1). To deal with the second, we use the method of [18; lemma 2] resting on the use of a Fejér kernel. We get, for a suitable constant $c_2 > 0$,

$$\mathfrak{R} \ll \frac{\Psi(x, y; f^\dagger)}{L_b(y)} + x\varrho_r(u)\mathrm{e}^{-c_2\overline{u}}.$$

Recall that $\mathcal{F}(s, y) := \mathcal{B}(s, y)\zeta(s, y)^z$. In order to bound the integral of (6·2), we need to estimate $\mathcal{F}(s, y)/\zeta(\alpha_r, y)^r$ with $z = r\mathrm{e}^{i\mathfrak{t}}$. Since $\log\zeta(2\alpha_r, y) \asymp \sqrt{\overline{u}}$, hypothesis (1·34) yields

$$\Big|\frac{\mathcal{F}(\alpha_r + i\tau, y)}{\zeta(\alpha_r, y)^r}\Big| \ll \mathrm{e}^{c_3\sqrt{\overline{u}} - rW(\tau, y)}, \quad \text{with} \quad W(\tau, y) := \sum_{p\leqslant y}\frac{1 - \cos(\mathfrak{t} - \tau\log p)}{p^{\alpha_r}}. \tag{6·3}$$

The sum $W(\tau, y)$ may then be estimated as in [18; lemma 1]. For $|\tau|\log y \leqslant \frac{1}{2}|\mathfrak{t}|$, we have

$$W(\tau, y) \gg \frac{1 - \cos(\mathfrak{t}/2)}{\log y}\sum_{p\leqslant y}\frac{\log p}{p^{\alpha_r}} \gg \overline{u}.$$

For the range $|\mathfrak{t}|/(2\log y) < |\tau| \leqslant L_{b_1}(y)$, we appeal to a variant of [20; lemma III.5.16] proved by contour integration using the Korobov-Vinogradov zero-free region of the Riemann zeta function. We get

$$W(\tau, y) \geqslant \sum_{p\leqslant y}\frac{\{1 - \cos(\mathfrak{t} - \tau\log p)\}\log p}{p^{\alpha_r}\log y} \gg \frac{\overline{u}\tau^2}{(1-\alpha_r)^2 + \tau^2} \gg \frac{\overline{u}}{(\log\overline{u})^2}.$$

Carrying back into (6·3), we obtain

$$\Big|\frac{\mathcal{F}(\alpha_r + i\tau, y)}{\zeta(\alpha_r, y)^r}\Big| \ll \mathrm{e}^{-c_4\overline{u}/(\log\overline{u})^2}.$$

The proof is now completed by inserting these bounds into the main term of Perron's formula (6·2), taking the lower bound (6·1) into account. The assertion regarding local uniformity follows from the fact that all bounds used in this proof are uniform in any compact subset of $H^- \smallsetminus \{0\}$.

## 7. Proofs of Corollaries 2.1, 2.2 and 2.3

Recall the definitions of $\sigma_r := \sigma_r(x,y)$ and $\mu_r := \mu_{x,y}$ in (2·1).

### 7·1. Preliminary estimates

Let $r$ be a positive, bounded parameter and $z := r\mathrm{e}^{it}$ $(t\in\mathbb{R})$. Put

$$u_t := u/z = \mathrm{e}^{-it}u/r,\quad w_t := \zeta_0(u_t),\quad F_r(t) := uw_t - r\mathrm{e}^{it}I(w_t),$$

so that, by (1·13) and (1·17),

$$\varrho_z(u) = \frac{\{1+O(1/u)\}\mathrm{e}^{\gamma z - F_r(t)}}{\sqrt{2\pi I''(w_t)}} = \varrho_r(u)\mathrm{e}^{F_r(0)-F_r(t)}\Big\{1+O\Big(|t|+\frac{1}{u}\Big)\Big\}, \tag{7·1}$$

$$R_z(u) \asymp \frac{\mathrm{e}^{-F_r(t)}}{\sqrt{u}}\cdot \tag{7·2}$$

For the sake of further reference, we note that

$$I(w_t) = u_t\Big\{1+\frac{1}{w_t}+\frac{2}{w_t^2}+O\Big(\frac{1}{\xi(u)^3}\Big)\Big\},\quad I'(w_t)=u_t,\quad I''(w_t) = u_t-(u_t-1)/w_t, \tag{7·3}$$

where the first estimate follows from (3·7).

Applying Corollary 1.3 to $f_z(n) := z^{\omega(n)}$ $(n\geqslant 1)$, we get, for $|t|<\pi/2$,

$$\begin{aligned}\Psi(x,y;f_z) &= a_0(f_z)x\varrho_z(u)(\log y)^{z-1} + O\Big(\frac{xR_z(u)\log 2u}{(\log y)^{2-z}} + \frac{x\mathbf{1}_{[1,1+e_y]}(u)}{\log y}\Big)\\ &= a_0(f_z)x\varrho_z(u)(\log y)^{z-1} + O\Big(\frac{xR_z(u)\log 2u}{\log y}\Big),\end{aligned} \tag{7·4}$$

with

$$a_0(f_z) = B(z) := \prod_p \Big(1-\frac{1}{p}\Big)^z\Big(1+\frac{z}{p-1}\Big).$$

We shall need extra precision in (7·1). Extending the proof of [17; th. 1], relevant to real $z$, to the case $\Re e\, z>0$, we get, for any fixed $J\geqslant 1$,

$$\varrho_z(u) = \frac{\mathrm{e}^{\gamma z-F_r(t)}}{\sqrt{2\pi I''(w_t)}}\Big\{1+\sum_{1\leqslant j\leqslant J}\frac{\lambda_j(w_t)}{I''(w_t)^j} + O\Big(\frac{1}{u^{J+1}}\Big)\Big\}$$

where the $\lambda_j$ $(1\leqslant j\leqslant J)$ are analytic functions bounded in a neighbourhood of $\mathbb{R}^+$ containing $w_t$. Bearing in mind that $w_0 = \xi(u/r)$, it follows that

$$\begin{aligned}\varrho_z(u) &= \frac{\{1+O(t)\}\mathrm{e}^{\gamma r-F_r(t)}}{\sqrt{2\pi I''(w_0)}}\Big\{1+\sum_{1\leqslant j\leqslant J}\frac{\lambda_j(w_0)}{I''(w_0)^j}+O\Big(\frac{|t|}{u}+\frac{1}{u^{J+1}}\Big)\Big\}\\ &= \varrho_r(u)\mathrm{e}^{F_r(0)-F_r(t)}\Big\{1+O\Big(|t|+\frac{1}{u^{J+1}}\Big)\Big\}.\end{aligned} \tag{7·5}$$

Next, let us perform the Taylor expansion of $F_r(t)+r\mathrm{e}^{it}\log_2 y$ to the second order at the origin. By [14; lemma 1] we have

$$\zeta_0'(v) = \frac{\zeta_0(v)}{v\{\zeta_0(v)-1\}}\Big\{1+O\Big(\frac{1}{v\xi(v)}\Big)\Big\}, \tag{7·6}$$

from which we derive

$$\begin{aligned} w_t' &= -iu_t\zeta_0'(u_t) = \frac{-iw_t}{w_t-1}\Big\{1+O\Big(\frac{1}{u\xi(u)}\Big)\Big\},\\ w_t'' &= -u_t\zeta_0'(u_t) - u_t^2\zeta_0''(u_t) = \frac{1}{w_t^2}+O\Big(\frac{1}{\xi(u)^3}\Big),\end{aligned} \tag{7·7}$$

where the second estimate is derived by differentiating [14; (4.2)] or by differentiating the definition $\mathrm{e}^{w_t} = 1+u_tw_t$ and using the first formula.

Now, by (7·3) and (7·7), we have, for large $u$,

$$
\begin{aligned}
(7{\cdot}8)\qquad F_r'(t) &= -ir\mathrm{e}^{it}I(w_t) = -iu\Big\{1+\frac{1}{w_t}+\frac{2}{w_t^2}+O\Big(\frac{1}{\xi(u)^3}\Big)\Big\},\\
F_r''(t) &= r\mathrm{e}^{it}I(w_t) - iuw_t' = r\mathrm{e}^{it}I(w_t) - uu_t\zeta_0'(u_t)\\
(7{\cdot}9)\qquad &= u\Big\{1+\frac{1}{w_t}+\frac{2}{w_t^2}+O\Big(\frac{1}{\xi(u)^3}\Big)\Big\} - \frac{uw_t}{w_t-1}\Big\{1+O\Big(\frac{1}{uw_0}\Big)\Big\}\\
&= \frac{u}{w_t^2}\Big\{1+O\Big(\frac{1}{w_0}\Big)\Big\},\\
(7{\cdot}10)\qquad F_r'''(t) &= ir\mathrm{e}^{it}I(w_t) + r\mathrm{e}^{it}w_t'I'(w_t) - iuw_t'' = \frac{iu}{w_t^2}\Big\{1+O\Big(\frac{1}{w_0}\Big)\Big\} - iuw_t'' \ll \frac{u}{w_0^3}\cdot
\end{aligned}
$$

It follows that

$$
(7{\cdot}11)\qquad
\begin{aligned}
H_r(t) &:= r(\mathrm{e}^{it}-1)\log_2 y + F_r(0) - F_r(t)\\
&= r(\mathrm{e}^{it}-1)\log_2 y - tF_r'(0) - \tfrac12 t^2F_r''(0) + O\big(\sigma^2t^3\big)\\
&= it\mu_r - \tfrac12 t^2\sigma_r^2 + O\big(t^3\sigma^2\big).
\end{aligned}
$$

Next, we need to estimate $\varrho_z(u)-\varrho_r(u)$ for small $t$.

As a start, we observe that it follows from the first equation in (7·7) that

$$
\frac{\mathrm{d}I(w_t)}{\mathrm{d}t} \ll u,
$$

whence $R_z(u)\ll\varrho_r(u)$ for $tu\ll 1$ and $u\geqslant 1$.

Now, under the assumption that $\Re e\, z>0$, we also have

$$
\begin{aligned}
\varrho_z(u) &= \frac{1}{2\pi i}\int_{1+i\mathbb{R}} \mathrm{e}^{\gamma z+zI(-s)+us}\,\mathrm{d}s,\\
\frac{\mathrm{d}\varrho_z(u)}{\mathrm{d}z} &= \frac{1}{2\pi i}\int_{1+i\mathbb{R}} \big\{\gamma+I(-s)\big\}\mathrm{e}^{\gamma z+zI(-s)+us}\,\mathrm{d}s\\
&= \gamma\varrho_z(u) + \frac{1}{2\pi i}\int_{1+i\mathbb{R}} \mathrm{e}^{\gamma z+zI(-s)+us}\int_0^1 \frac{\mathrm{e}^{-hs}-1}{h}\,\mathrm{d}h\,\mathrm{d}s\\
&= \gamma\varrho_z(u) + \frac{1}{2\pi i}\int_0^1\frac{\mathrm{d}h}{h}\int_{1+i\mathbb{R}} (\mathrm{e}^{-hs}-1)\mathrm{e}^{\gamma z+zI(-s)+us}\,\mathrm{d}s\\
&= \gamma\varrho_z(u) + \int_0^1 \frac{\varrho_z(u-h)-\varrho_z(u)}{h}\mathrm{d}h \ll \int_0^1 \varrho_z(u)\xi(u)\mathrm{e}^{h\xi(u/r)}\,\mathrm{d}h \ll u\varrho_r(u).
\end{aligned}
$$

We therefore obtain

$$
(7{\cdot}12)\qquad \varrho_z(u)-\varrho_r(u)\ll tu\varrho_r(u)\qquad (z=r\mathrm{e}^{it},\ tu\ll 1).
$$

### *7·2. Proof of Corollary 2.1*

It follows from (7·4) and (7·5) with $r=1$ that

$$
(7{\cdot}13)\qquad \Psi(x,y;f_z) = \Big\{1+O\Big(|t|+\frac{1}{u^{J+1}}\Big)\Big\}\Psi(x,y)\mathrm{e}^{-H_1(t)} + O\Big(\frac{\Psi(x,y)\log 2u}{\log y}\Big)
$$

and hence, for $|t|\leqslant\sigma^{2/3}$,

$$
(7{\cdot}14)\qquad
\begin{aligned}
h_{x,y}(t) &:= \frac{1}{\Psi(x,y)}\sum_{n\in S(x,y)} \mathrm{e}^{it\{\omega(n)-\mu\}/\sigma}\\
&= \Big\{1+O\Big(\frac{1}{u^{J+1}}+\frac{|t|(1+t^2)}{\sigma}\Big)\Big\}\mathrm{e}^{-t^2/2} + O\Big(\frac{\log 2u}{\log y}\Big).
\end{aligned}
$$

This will be of use for intermediate values of $t$. For large values of $|t|$, we deduce from the above that, with a suitable constant $c > 0$, we have

$$h_{x,y}(t) \ll \mathrm{e}^{-t^2/3} \qquad (|t| \leqslant c\sigma). \tag{7·15}$$

For small values of $t$, we appeal to the estimate

$$\frac{1}{\Psi(x,y)} \sum_{n\in S(x,y)} \{\omega(n)-\mu\}^2 \ll \mu$$

resulting from theorem 1.1, lemma 9.1, and formula (9.4) of [4]. This implies

$$h_{x,y}(t) = 1 + O\Big(\frac{|t|\sqrt{\mu}}{\sigma} + \frac{t^2\mu}{\sigma^2}\Big) = \mathrm{e}^{-t^2/2} + O\Big(\frac{|t|\sqrt{\mu}}{\sigma} + \frac{t^2\mu}{\sigma^2}\Big) \qquad (t \ll 1). \tag{7·16}$$

By the Berry-Esseen inequality (see, e.g., [20; th. II.7.16]), we therefore obtain

$$\frac{1}{\Psi(x,y)} \sum_{\substack{n\in S(x,y)\\ \omega(n)-\mu\leqslant v\sigma}} 1 = \Phi(v) + O(R) \qquad (v\in\mathbb{R}), \tag{7·17}$$

with

$$R := \frac{1}{T} + \int_{-T}^{T} \big|h_{x,y}(t) - \mathrm{e}^{-t^2/2}\big| \frac{\mathrm{d}t}{t}\cdot$$

Selecting $T := \min(\sigma^2, u^{J+1})$, we hence get

$$\begin{aligned} R &\ll \int_0^{1/T^2} \Big(\frac{\sqrt{\mu}}{\sigma} + \frac{t\mu}{\sigma^2}\Big)\,\mathrm{d}t + \int_{1/T^2}^{T^{2/3}} \Big(\frac{1}{tu^{J+1}} + \frac{1+t^2}{\sigma}\Big)\frac{\mathrm{d}t}{\mathrm{e}^{t^2/2}} + \int_{T^{2/3}}^{cT} \frac{\mathrm{d}t}{t\mathrm{e}^{t^2/3}} + \frac{1}{T} \\ &\ll \frac{\sqrt{\mu}}{T^2\sigma} + \frac{\log T}{u^{J+1}} + \frac{1}{\sigma} + \frac{1}{T} \ll \frac{\log 2u}{u^{J+1}} + \frac{1}{\sigma}, \end{aligned} \tag{7·18}$$

since $\sqrt{\mu}/\sigma \ll \log 2u$ and $\sigma^2 \leqslant u \Rightarrow \mu \ll \sigma^2$.

At this stage we haven't utilised (7·12). For $tu \ll \sigma$ and $t^3 \ll \sigma$, we derive from this approximation and (7·4) that, with now $z := \mathrm{e}^{it/\sigma}$,

$$\begin{aligned} &\sum_{n\in S(x,y)} \mathrm{e}^{it\{\omega(n)-\mu\}/\sigma} \\ &\qquad = x\varrho(u)\Big\{1 + O\Big(\frac{tu}{\sigma}\Big)\Big\}\mathrm{e}^{it(\log_2 y-\mu)/\sigma - t^2(\log_2 y)/2\sigma^2 + O(t^3/\sigma)} + O\Big(\frac{R_z(u)\log 2u}{\log y}\Big) \\ &\qquad = \Psi(x,y)\Big\{1 + O\Big(\frac{\log 2u}{\log y} + \frac{|t|u}{\sigma}\Big)\Big\}\mathrm{e}^{-t^2/2}\Big\{1 + O\Big(\frac{|t|u}{\sigma} + \frac{t^2 u}{(\sigma\log 2u)^2} + \frac{|t|^3}{\sigma}\Big)\Big\}, \end{aligned}$$

and hence

$$h_{x,y}(t) = \mathrm{e}^{-t^2/2}\Big\{1 + O\Big(\frac{\log 2u}{\log y} + \frac{|t|u}{\sigma} + \frac{t^2u}{(\sigma\log 2u)^2} + \frac{|t|^3}{\sigma}\Big)\Big\}.$$

We thus have

$$\int_{1/\sigma}^{\sigma} \Big|h_{x,y}(t) - \mathrm{e}^{-t^2/2}\Big| \frac{\mathrm{d}t}{t} \ll \frac{(\log X)\log 2u}{\log y} + \frac{u}{\sigma} \ll \frac{u}{\sigma}\cdot$$

Appealing to (7·16) in order to bound the corresponding integral over the range $0 \leqslant t \leqslant 1/\sigma$, we see that the remainder $R$ arising in (7·17) satisfies

$$R \ll \frac{u}{\sigma}\cdot \tag{7·19}$$

Gathering our two bounds for $R$, we derive (2·2).

### 7·3. Proof of Corollary 2.2

Let $r>0$, $z:=r\mathrm{e}^{it}$, $f_z(n):=z^{\omega(n)}$, and

$$j_{x,y}(z):=\frac{\Psi(x,y;f_z)}{r^k\Psi(x,y)}$$

so that, with notation (2·3), we have

$$\psi_k(x,y)=\frac{1}{2\pi}\int_{-\pi}^{\pi} j_{x,y}(r\mathrm{e}^{it})\mathrm{e}^{-ikt}\,\mathrm{d}t.$$

Recall (7·4). Carrying (7·11) back into (7·5) and (7·4), we get, for $|t|\leqslant \sigma_r^{-2/3}$,

$$j_{x,y}(r\mathrm{e}^{it})\mathrm{e}^{-itk}=A_{u,y}(r,k)\bigg(\Big\{1+O\big(E_{u,y}(t)\big)\Big\}\mathrm{e}^{it(\mu_r-k)-\frac12 t^2\sigma_r^2}+O\Big(\frac{\log 2u}{\log y}\Big)\bigg),\tag{7·20}$$

with

$$A_{u,y}(r,k):=\frac{B(r)\varrho_r(u)(\log y)^{r-1}}{r^k\varrho(u)},\quad E_{u,y}(t):=\frac{1}{u^{J+1}}+|t|+|t|^3\sigma^2.$$

At this stage, we deploy two distinct strategies according to whether we want to prove statement (a) or (b).

To establish (a), we select $r=1$ and integrate over $t\in[-\sigma^{-2/3},\sigma^{-2/3}]$. We get

$$\frac{1}{2\pi}\int_{-1/\sigma^{2/3}}^{1/\sigma^{2/3}} j_{x,y}(\mathrm{e}^{it})\mathrm{e}^{-ikt}\,\mathrm{d}t=\frac{\mathrm{e}^{-\frac12(k-\mu)^2/\sigma^2}}{\sqrt{2\pi}\sigma}+\mathcal{R}_1$$

with

$$\mathcal{R}_1\ll\frac{1}{u^{J+1}\sigma}+\frac{1}{\sigma^2}\cdot\tag{7·21}$$

To improve precision in the case of small values of $u$, we employ (7·12). When $tu\ll 1$, the main term of (7·4) is, writing $\beta(t):=it(\log_2 y-k)-\frac12 t^2\log_2 y$,

$$x\varrho(u)\big\{1+O\big(tu\big)\big\}\mathrm{e}^{\beta(t)+O(t^3\log_2 y)}=\Psi(x,y)\Big\{1+O\Big(\frac{\log 2u}{\log y}+|t|u\Big)\Big\}\mathrm{e}^{\beta(t)}\big\{1+O\big(t^3\sigma^2\big)\big\}.$$

Since $\mu-\log_2 y\ll u$, $\sigma^2-\log_2 y\ll u$, we have $\beta(t)=i\mu t-\frac12 t^2\sigma^2+O(|t|u+t^2u)$, and so we get

$$j_{x,y}(\mathrm{e}^{it})\mathrm{e}^{-ikt}=\mathrm{e}^{it(\mu-k)-\frac12 t^2\sigma^2}\bigg\{1+O\bigg(\frac{\log 2u}{\log y}+|t|u+t^2u+|t|^3\sigma^2\bigg)\bigg\}.$$

Under the extra assumption $u\leqslant\sigma^{2/3}$, we therefore obtain

$$\int_{-1/\sigma^{2/3}}^{1/\sigma^{2/3}} j_{x,y}(\mathrm{e}^{it})\mathrm{e}^{-ikt}\,\mathrm{d}t=\frac{\mathrm{e}^{-\frac12(k-\mu)^2}}{\sqrt{2\pi}\sigma}+O(\mathcal{R}_1),$$

with

$$\mathcal{R}_1\ll\frac{u}{\sigma^2}\cdot\tag{7·22}$$

Gathering (7·21) and (7·22), we obtain

$$\mathcal{R}_1\ll\frac{u}{\sigma^2+u^{J+2}\sigma}+\frac{1}{\sigma^2}\cdot$$

In order to complete the proof of (2·5), we need to estimate

$$E:=\int_{1/\sigma^{2/3}<|t|\leqslant\pi} j_{x,y}(\mathrm{e}^{it})\mathrm{e}^{-ikt}\,\mathrm{d}t.$$

We shall obtain the required result from an upper bound for $\Psi(x,y;f_z)/\Psi(x,y)$.

Let $c\in]0,\frac12\pi[$. If $\frac12\pi-c<|t|\leqslant\pi$, we have $\Re e\, z\leqslant\sin c$. Taking account of the estimates (3·12) and (3·16), the upper bound derived from Corollary 1.3 provides

$$\Psi(x,y;f_z)/\Psi(x,y)\ll 1/(\log y)^{1-\sin c}$$

which is more than sufficient.

When $1/\sigma^{1/3}<|t|\leqslant\frac12\pi-c$, we have

$$\frac{\Psi(x,y;f_z)}{\Psi(x,y)}\ll \mathrm{e}^{-H(t)}$$

and the estimate (7·9) for $F''(t)$ implies

$$\begin{aligned}\Re e\, H(t)&=-(1-\cos t)\log_2 y-\int_0^t\Re e\, F''(v)(t-v)\,\mathrm{d}v\\&\leqslant-\frac{2t^2}{\pi^2}\Big(\log_2 y+\frac{2u}{\xi(u)^2}\Big)+O(1)\leqslant-\frac{2t^2\sigma^2}{\pi^2}+O(1).\end{aligned}$$

We therefore infer in turn

$$\frac{\Psi(x,y;f_z)}{\Psi(x,y)}\ll\mathrm{e}^{-2t^2\sigma^2/\pi^2},\qquad E\ll\frac{1}{\sigma^2}\cdot$$

This finishes the proof of (2·5).

Let us now embark on the proof of statement (b) and so select $r$ in (7·20) such that $\mu_r=k$, i.e. $r=k/L$ with notation (2·4). By (3·7), this is plainly possible with bounded $r$ under condition (2·6). We infer that

$$\frac{1}{2\pi}\int_{-1/\sigma_r^{2/3}}^{1/\sigma_r^{2/3}}j_{x,y}(r\mathrm{e}^{it})\,\mathrm{d}t=\frac{A_{u,y}(r,k)}{\sqrt{2\pi}\sigma_r}\big\{1+O(\mathcal{R}_2)\big\}\qquad(v\in\mathbb{R}),$$

with

$$\mathcal{R}_2\ll\frac{1}{u^{J+1}}+\frac{1}{\sigma}\cdot\tag{7·23}$$

Moreover, for this choice of $r$, we have, by Stirling's formula,

$$\frac{A_{u,y}(r,k)}{\sqrt{2\pi}\sigma_r}=\frac{B(r)\varrho_r(u)L^k\mathrm{e}^k\mathrm{e}^{-rI(\xi(u/r))}\sqrt{\mu_r}}{\sigma_r\varrho(u)k^k(\log y)\sqrt{2\pi k}}=\mathcal{K}_r(u)\frac{\mathrm{e}^{-L}L^k}{k!}\Big\{1+O\Big(\frac1\mu\Big)\Big\}$$

To improve precision in the case of small values of $u$, we again appeal to (7·12). When $tu\ll1$, the main term of (7·4) is, writing $\beta_r(t):=it(r\log_2 y-k)-\frac12t^2\log_2 y$,

$$\begin{aligned}&xA_{u,y}(r,k)\varrho(u)\big\{1+O\big(tu\big)\big\}\mathrm{e}^{\beta_r(t)+O(t^3\log_2 y)}\\&\qquad=\Psi(x,y)A_{u,y}(r,k)\Big\{1+O\Big(\frac{\log 2u}{\log y}+|t|u\Big)\Big\}\mathrm{e}^{\beta(t)}\big\{1+O\big(t^3\sigma^2\big)\big\}.\end{aligned}$$

Since $\mu_r-\log_2 y\ll u$, $\sigma_r^2-\log_2 y\ll u$, we have $\beta_r(t)=it\mu_r-\frac12t^2\sigma_r^2+O\big(|t|u+t^2u\big)$, and so we get

$$j_{x,y}(r\mathrm{e}^{it})\mathrm{e}^{-ikt}=A_{u,y}(r,k)\mathrm{e}^{-\frac12t^2\sigma_r^2}\Big\{1+O\Big(\frac{\log 2u}{\log y}+|t|u+t^2u+|t|^3\sigma^2\Big)\Big\}.$$

Under the extra assumption $u\leqslant\sigma_r^{2/3}$, we therefore obtain, still with $\mu_r=k$,

$$\int_{-1/\sigma_r^{2/3}}^{1/\sigma_r^{2/3}}j_{x,y}(r\mathrm{e}^{it})\mathrm{e}^{-ikt}\,\mathrm{d}t=\frac{A_{u,y}(r,k)}{\sqrt{2\pi}\sigma_r}\big\{1+O(\mathcal{R}_2)\big\},$$

with

$$\mathcal{R}_2\ll\frac{u}{\sigma}\cdot\tag{7·24}$$

Gathering (7·23) and (7·24), we obtain

$$\mathcal{R}_2\ll\frac{u}{\sigma+u^{J+2}}+\frac1\sigma\cdot$$

As previously, we complete the proof of (2·7) by bounding

$$E_r:=\int_{1/\sigma_r^{2/3}<|t|\leqslant\pi}j_{x,y}(r\mathrm{e}^{it})\mathrm{e}^{-ikt}\,\mathrm{d}t.$$

The required estimate follows from the upper bounds for $\Psi(x,y;f_z)/\Psi(x,y)$ derived from Corollary 1.3 in much the same manner than in the proof of (2·5) and we omit the details. This finishes the proof of (2·7).

### *7·4. Proof of Corollary 2.3*

Let $r \in [\frac{1}{2}, 2]$. It follows from (1·13) that $\log\{\varrho_r(u)/\varrho(u)\} = G(r) + O(1)$ with

$$G(r) := u\xi(u) - u\xi(u/r) + rI(\xi(u/r)) - I(\xi(u)),$$

so that

$$G(1) = 0, \quad G'(r) = I(\xi(u/r)), \quad G''(r) = -\xi'(u/r)u^2/r^3, \qquad G'''(r) \ll u.$$

For $|h| \leqslant \frac{1}{2}$, we therefore have

$$G(1+h) = hI(\xi(u)) - \tfrac{1}{2}\xi'(u)u^2h^2 + O\big(uh^3\big).$$

Now, for $\frac{1}{2} \leqslant r \leqslant 1$, $v \geqslant 0$, we have

$$\frac{1}{\Psi(x,y)} \sum_{\substack{n \in S(x,y)\\ \omega(n) \leqslant \mu - v\sigma}} 1 \leqslant \frac{r^{-\mu+v\sigma}}{\Psi(x,y)} \sum_{n \in S(x,y)} r^{\omega(n)} \ll r^{-\mu+v\sigma}(\log y)^{r-1}\mathrm{e}^{G(r)}.$$

Writing $r = 1 - h$, $0 \leqslant h \leqslant \frac{1}{2}$, we derive that

$$\begin{aligned}(v\sigma - \mu)\log r - h\log_2 y + G(r)&\\ &= h\{\mu - v\sigma - \log_2 y - I(\xi)\} - \tfrac{1}{2}h^2\{v\sigma - \mu + u^2\xi'(u)\} + O\big(h^3u\big)\\ &= -hv\sigma - \tfrac{1}{2}h^2(v\sigma - \sigma^2) + O(h^3u)\end{aligned}$$

Selecting $h = v/\sigma$, we get

$$\frac{1}{\Psi(x,y)} \sum_{\substack{n \in S(x,y)\\ \omega(n) \leqslant \mu - v\sigma}} 1 \ll \mathrm{e}^{-\frac{1}{2}v^2 + O(v^3u/\sigma^3)} \quad (v \geqslant 0).$$

This bound is certainly non-trivial if $v \leqslant c_1\sigma/(\log 2u)^2$ with suitable absolute $c_1 > 0$.

Arguing similarly with $r = 1 + h \geqslant 1$, we get

$$\frac{1}{\Psi(x,y)} \sum_{\substack{n \in S(x,y)\\ \omega(n) \geqslant \mu + v\sigma}} 1 \ll \mathrm{e}^{-\frac{1}{2}v^2 + O(uv^3/\sigma^3 + v^3/\sigma)} \qquad (v \geqslant 0).$$

The two previous estimates imply (2·8): for sufficiently small $c_1$ and $v = v_0 := c_1\sigma/(\log 2u)^2$, the bounds are $\ll \mathrm{e}^{-\frac{1}{3}c_1^2\sigma^2/(\log 2u)^4}$, for $v \leqslant v_0$, they are $\ll \mathrm{e}^{-v^2/3}$. If, furthermore, we have $v \ll \min\{\sigma/u^{1/3}, \sigma^{1/3}\}$, they are $\ll \mathrm{e}^{-v^2/2}$.

**Acknowledgement.** The authors take pleasure in addressing warm thanks to the referee whose thorough remarks and suggestions led to significant improvements of this work.

Régis de la Bretèche
Université Paris Cité, Sorbonne Université, CNRS
Institut de Math. de Jussieu-Paris Rive Gauche
F-75013 Paris
France
regis.de-la-breteche@imj-prg.fr

Gérald Tenenbaum
Institut Élie Cartan
Université de Lorraine
BP 70239
54506 Vandœuvre-lès-Nancy Cedex
France
gerald.tenenbaum@univ-lorraine.fr